\documentclass{amsart}
\usepackage{color}
\usepackage{graphicx,tikz}
\usepackage{multicol}

\newcommand{\R}{\ensuremath{\mathbb{R}}}
\newcommand{\Z}{\ensuremath{\mathbb{Z}}}

\renewcommand{\H}{\ensuremath{\mathbb{H}}}

\newcommand{\Ccell}{C^{\text{cell}}}
\newcommand{\CLip}{C^{\text{Lip}}}
\newcommand{\Zcell}{Z^{\text{cell}}}
\newcommand{\ZLip}{Z^{\text{Lip}}}
\renewcommand{\div}{\mathop{\rm Div}\nolimits}
\newcommand{\ddiv}{\mathop{\rm div}\nolimits}

\newcommand{\Lip}{\mathop{\rm Lip}}

\newcommand{\abs}{\mathop{\sf abs}}

\newcommand{\pos}{\mathop{\sf pos}}
\renewcommand{\neg}{\mathop{\sf neg}}
\newcommand{\id}{\mathop{id}}
\newcommand{\Orth}{\mathop{\sf Orth}\nolimits}
\newcommand{\supp}{\mathop{\rm supp}}

\renewcommand{\P}{\ensuremath{\mathcal{P}}}
\newcommand{\QQ}{\ensuremath{\mathcal{Q}}}

\newcommand{\mass}{\mathop{\rm mass}}
\newcommand{\V}{l}
\newcommand{\divdim}{\mathop{\rm divdim}}

\newcommand{\Agamma}{\ensuremath{A_\Gamma}}
\newcommand{\Xgamma}{\ensuremath{X_\Gamma}}
\newcommand{\Hgamma}{\ensuremath{H_\Gamma}}

\theoremstyle{plain}
\newtheorem{theorem}{Theorem}[section]
\newtheorem{thm}[theorem]{Theorem}
\newtheorem{lemma}[theorem]{Lemma}
\newtheorem{proposition}[theorem]{Proposition}
\newtheorem{prop}[theorem]{Proposition}

\numberwithin{equation}{section}

\theoremstyle{definition}
\newtheorem{definition}[theorem]{Definition}

\newtheorem{remark}[theorem]{Remark}

\begin{document}

\title{Pushing fillings in  right-angled Artin groups}
\author[A. Abrams]{Aaron Abrams}
\address{
Aaron Abrams\\
Dept.\ of Math and Comp. Sci.\\
Emory University\\
Atlanta, GA 30322\\ 
USA}
\email{abrams@mathcs.emory.edu}
\author[N. Brady]{Noel Brady}
\address{
Noel Brady\\
Department of Mathematics\\
University of Oklahoma\\
Norman, OK 73019\\ 
USA}
\email{nbrady@math.ou.edu}
\author[P. Dani]{Pallavi Dani}
\address{
Pallavi Dani\\
Department of Mathematics\\
Louisiana State University\\
Baton Rouge, LA 70803\\
USA}
\email{pdani@math.lsu.edu}
\author[M. Duchin]{Moon Duchin}
\address{
Moon Duchin\\
Department of Mathematics\\
University of Michigan\\
Ann Arbor, MI 48109\\
USA}
\email{mduchin@umich.edu}
\author[R. Young]{Robert Young}
\address{
Robert Young\\
Department of Mathematics\\
University of Toronto\\
Toronto, ON  M5S 2E4\\
Canada
}
\email{rjyoung1729@gmail.com}
\date{\today}
\thanks{This work was supported by a {\sf SQuaRE} grant from the
  American Institute of Mathematics.  The second author and 
  the fourth author are partially supported by NSF grants DMS-0906962 and 
  DMS-0906086, respectively.
  The fifth author would like to
  thank New York University for its hospitality during the preparation
of this paper.}

\begin{abstract}
We construct {\em pushing maps} on the cube complexes that model
right-angled Artin groups (RAAGs) in order to study  filling
problems in certain subsets of these cube complexes.  We use
radial pushing  to obtain upper bounds on {\em higher
divergence functions}, finding that the $k$-dimensional divergence of
a RAAG is bounded by $r^{2k+2}$.  
These divergence functions, previously defined for Hadamard manifolds
to measure isoperimetric properties ``at infinity,"
are defined here as a family of quasi-isometry invariants of groups.
By pushing along the height gradient, we
also show
that the $k$-th order Dehn function of a Bestvina-Brady group is
 bounded by $V^{(2k+2)/k}$.  We
construct a class of RAAGs called {\em orthoplex groups} which show
that each of these upper bounds is sharp.
\end{abstract}

\maketitle

\section{Introduction}

Many of the groups studied in geometric group theory are subgroups of
non-positively curved groups. This family includes lattices in
symmetric spaces, Bestvina-Brady groups, and many solvable groups. In
each of these cases, the group acts geometrically on a subset of a non-positively
curved space, and one can approach the study of the subgroup by
considering the geometry of that part of the model space.

In this paper, we construct {\em pushing maps} for the cube complexes
that model right-angled Artin groups.  These maps serve to modify chains
so that they lie in special subsets of the space.  
We find that the geometry of the groups is ``flexible" enough that 
it is not much more difficult to fill curves and cycles
in these special subsets than to fill them efficiently in the ambient
space.  One application will be to study higher divergence functions,
which measure the geometry of a group ``at infinity" by avoiding a
large ball around the origin.  Another will be to push along the height
gradient in order to solve higher-order filling problems in Bestvina-Brady groups.

{\em Right-angled Artin groups} (or RAAGs) are given by presentations in which
each relator is a commutator of two generators.  A great deal is known about the
algebra, geometry, and combinatorics of RAAGs.  For instance, they have automatic 
structures and useful 
normal forms, and they act geometrically on CAT(0) cube complexes.  Many
tools are available for their study, in part because these
complexes contain flats arising from mutually commuting
elements  (see  \cite{charney}).  
Frequently, topological invariants of RAAGs can be read off 
of the defining graph, and along these lines we will relate properties of the
graph to the filling functions and divergence functions of the groups.

Filling functions describe the difficulty of finding a disc or chain
with a given boundary.  Recall that
the most basic filling function in groups is the {\em Dehn function}, which
measures the area necessary to fill a closed loop by a disk; 
these functions have been a key part of geometric group theory
since Gromov used them to characterize hyperbolic groups
(or arguably longer, since Dehn used related ideas to find fast solutions
to the word problem).
This can naturally be generalized to higher-order
Dehn functions, which describe the difficulty of filling $k$-spheres
by $(k+1)$-balls or $k$-cycles by $(k+1)$-chains.

{\em Topology at infinity} is the study of the asymptotic structure of
groups by attaching topological invariants to the complements of large
balls; the theorems of Hopf and Stallings 
about ends of groups were early examples of major results of this kind.

One can study this topology at infinity quantitatively by
introducing filling invariants at infinity, such as the {\em higher divergence
functions}, which measure rates of filling in complements of large balls in groups
and other metric spaces.  With respect to a fixed basepoint $x_0$ in a
space, we will describe a map whose image is disjoint from the ball of
radius $r$ about $x_0$ as being $r$-{\em avoidant}.  Roughly, the
$k$-dimensional divergence function is a filling invariant for
avoidant cycles and chains (or spheres and balls); it measures the volume 
needed to fill an avoidant $k$-cycle by an avoidant $(k+1)$-chain.  
(We will make this precise in \S\ref{sec:defineDivergence}.)
As with Dehn functions, the divergence functions become meaningful for
finitely generated groups by adding an appropriate equivalence
relation to make the definition
invariant under quasi-isometry.

For $k\ge 0$, our functions $\div^k$ are closely related to the higher 
divergence functions defined by 
Brady and Farb in \cite{brady-farb} for the special case of Hadamard manifolds.
Using the manifold definition, combined results of Leuzinger and Hindawi prove
that the higher divergence functions detect the real-rank of 
a symmetric space, as Brady-Farb  had conjectured \cite{leuzinger, hindawi}.  
Thus the geometry and the algebra are connected.
Wenger generalized this, showing that higher divergence functions
detect the Euclidean rank of any CAT(0) space \cite{wenger2}.  In 
extending these notions to groups, this paper is necessarily largely concerned with precise definitions
and with tractable cases, but it may be regarded as making the first steps in a process of discovering
which properties of groups are detected by this family of invariants.  A secondary contribution of the present
paper is in providing  in \S\ref{sec:compare}
what we hope is  a brief but usable treatment of the comparison between 
the various categories of 
filling functions found in the literature; here, we 
explain why the main results and techniques in this paper,
though their properties are stated and proved in the homological category,
work just as well with homotopical definitions.

Higher divergence functions are interesting in part because they unify a number of concepts in 
coarse geometry and geometric group theory.  For instance, in the $k=0$ case, this is the classical
{\em divergence of geodesics}, which relates to the curvature and in particular detects  hyperbolicity.  
(Gromov showed that a space is $\delta$-hyperbolic iff it has 
exponential divergence of geodesics in a certain precise sense.)  
More recently, many papers in 
geometric group theory have studied polynomial divergence of geodesics, including but not limited to
\cite{gersten1, gersten2, kapovich-leeb, duchin-rafi, macura, d-m-s}.  Much of this work arose 
to explore an expectation offered by Gromov in \cite{groAII} that nonpositively curved
spaces should, like symmetric spaces, exhibit a gap in the possible rates of 
divergence of geodesics (between linear and exponential).  
On the contrary, it is now clear that quadratic
divergence of geodesics often occurs in groups where many ``chains of
flats'' are present,
and Macura has produced examples of groups with polynomial divergence
of every degree \cite{macura}.




In this paper we develop several applications of pushing maps (defined
in \S\ref{sec:pushing}), which are ``singular retractions'' defined from the complex 
associated to a RAAG onto various subsets of this complex, such as the
exterior of a ball or a Bestvina-Brady subgroup.  We will use these maps to
obtain upper bounds on higher Dehn functions of
Bestvina-Brady groups by pushing fillings into these subgroups (\S\ref{sec:BB}), and 
we will use them in a different way to find special examples 
called {\em orthoplex groups} where the upper bounds are achieved.
In \S\ref{sec:higherdiv} we
study higher divergence functions by pushing fillings
out of balls:  if a RAAG  $\Agamma$
is $k$-connected at infinity, we can guarantee that avoidant fillings satisfy a polynomial
bound, namely $\div^k(\Agamma) \preceq r^{2k+2}$.
Next, these upper bounds are shown to be sharp by using the earlier estimates
for orthoplex groups.

Although the upper and lower bounds are sharp in every dimension, we are not able to
 specify which rates of divergence occur between the two extremes
for $k\ge 1$.  However, for $k=0$, we show in \S\ref{sec:div0}
 that every RAAG must have either linear
or quadratic divergence, depending only on whether the group is a direct product
(a property that can be read off of the defining graph).

We note that sorting right-angled Artin groups by their 
``divergence spectra'' gives a tool for distinguishing quasi-isometry types;
the QI classification problem for RAAGs still has many outstanding cases, particularly
in higher dimensions.

Several of the techniques  developed to deal with RAAGs
have applications in  other groups.  
``Pushing'' may be used in the torsion analogs of Baumslag's metabelian group
to find that the Dehn function is at worst quartic, as shown in \cite{kass-riley}---a priori, 
it takes nontrivial work even to show that the Dehn function is polynomial.  
(In fact, it turns out to be quadratic, as shown by de Cornulier-Tessera in \cite{cor-tess}.)
These ideas are also applicable in so-called ``perturbed RAAGs,'' as explained in \cite{b-g-l}.
Finally, pushing techniques can be adapted to give results for divergence in mapping class 
groups, which we will explore in a future paper.


\section{Dehn functions and divergence functions}

In this section, we will define the higher-order Dehn functions and
the higher divergence functions and illustrate the basic methods of
this paper by using a pushing map to bound the divergence functions of
$\R^d$.

\subsection{Higher Dehn functions}\label{sec:compare}
We will primarily use homological Dehn functions,
following \cite{wordproc, groAII, Wenger}.  Homological Dehn functions
describe the difficulty of filling cycles by chains.  In contrast,
some other papers (\cite{AlWaPr, snowflake}) use homotopical Dehn
functions, which measure fillings of spheres by balls.

In low dimensions, these functions may differ, but they are
essentially the same for high-dimensional fillings in highly-connected
spaces.  If $X$ is a $k$-connected space and $k\ge 3$, the topologies of the boundary and of the
filling are irrelevant, and the homological and homotopical $k$-th order 
Dehn functions  of $X$ are the same.
When $k=2$, the topology of the boundary
is relevant, but the topology of the filling is not: a homological
filling of a sphere guarantees a homotopical one of nearly the same volume and vice versa, so that the homotopical
Dehn function is bounded above by the homological one \cite{groft1,groft2}.  (See also 
\cite[App.2.(A')]{groFRM}, \cite[Rem.2.6(4)]{snowflake}.)
The reverse is not true; spheres can be filled equally well by balls
or by chains, but there may exist cycles that are ``harder to fill''
than spheres \cite{YoungHomological}, and the homological
second-order Dehn function may be larger than the homotopical one.

The bounds in this paper on rates of filling of Lipschitz chains by Lipschitz cycles---for 
Euclidean space (Proposition~\ref{prop:eucBounds}), Bestvina-Brady groups 
(Theorems~\ref{thm:bestBradUpper},\ref{dehnsharp}), and right-angled Artin groups
(Theorem~\ref{divk} and the propositions used to prove it)---are all valid using homotopical definitions of the Dehn and divergence 
functions.
It is automatic that higher-order Dehn
function upper bounds stated for homological filling hold for
homotopical filling as well, for the general reasons given above.
An extra argument is needed 
in dimension 1 (for instance, to see that our upper bounds
on $\delta_G$ for Bestvina-Brady groups and on $\div^1$ for
right-angled Artin groups also hold in the homotopic category).  
Because our techniques below construct disks filling curves
rather than just chains, they also bound the homotopical Dehn
function (see also Remark~\ref{push-admiss}).  Likewise, the lower bounds that we prove use only spheres
as boundaries, so they hold equally well in both contexts.

We will define the higher Dehn function in two ways, one better-suited
to dealing with complexes, and one better for dealing with manifolds.

We define a polyhedral complex to be a CW-complex in which each cell
is isometric to a convex polyhedron in Euclidean space and the gluing
maps are isometries.  
If $X$ is a polyhedral complex, we can define filling
functions of $X$ based on cellular homology.  Assume that $X$
is $k$-connected and let $\Ccell_k(X)$ be the set of cellular
$k$-chains of $X$ with integer coefficients.  If $a\in \Ccell_k(X)$,
then $a=\sum_{i} a_i \sigma_i$ for some integers $a_i$ and distinct
$k$-cells $\sigma_i$.  Set
$\|a\|=\sum |a_i|$.
If $\Zcell_k(X)$ is the set of cellular $k$-cycles and $a\in
\Zcell_k(X)$, then the fact that $X$ is $k$-connected implies
that $a=\partial b$ for some $b\in \Ccell_{k+1}(X)$.  Define the
{\em filling volume}  and  the {\em $k$-th order Dehn function} 
to be
$$\delta^{k;\text{cell}}_{X}(a)=\mathop{\min_{b\in
    \Ccell_{k+1}(X)}}_{\partial b=a}\|b\|, \quad ; \qquad 
\delta^{k;\text{cell}}_{X}(\V) = \mathop{\max_{a\in
    \Zcell_{k+1}(X)}}_{\|a\|\le\V}\delta^{k;\text{cell}}_{X}(a).
$$

Alonso, Wang, and Pride \cite{AlWaPr} showed that if $G$ and $G'$ are
quasi-isometric groups acting geometrically (i.e., properly
discontinuously, cocompactly, and by isometries) on certain associated $k$-connected complexes $X$ and
$Y$ respectively, then $\delta^{k;\text{cell}}_X$ and
$\delta^{k;\text{cell}}_{Y}$ grow at the same rate; in particular,
this shows that the growth rate of $\delta^{k;\text{cell}}_X$
depends only on $G$, so we can define $\delta^{k;\text{cell}}_G$.  This is made rigorous by 
defining an equivalence relation $\asymp$ on functions, as follows.  There is a
partial order on the set of functions $\R^+\to \R^+$ given by the following symbol:
$f\preceq g$ means there exists $A>0$ such that
$$f(t) \le Ag(At+A)+At+A$$
for all $t\ge 0$, and the same property may be written $g \succeq f$.
Then $f\asymp g$ if and only
if $f\preceq g$ and $f \succeq g$.  This is the standard notion of
equivalence for coarse geometry, because it amounts to allowing a
linear rescaling of domain and range, as in a quasi-isometry.  Note
that the equivalence relation $\asymp$ identifies all linear and
sublinear functions into one class, but distinguishes polynomials of
different degrees.

Another way to define a homological higher-order Dehn function,
somewhat better suited to Riemannian manifolds and CAT(0)-spaces, is
to use singular Lipschitz chains, as in \cite{wordproc},
\cite{groFRM}, and \cite{Wenger}.  A full introduction to this approach
can be found in Chapter 10.3 of \cite{wordproc}.  Let $X$ be a
$k$-connected Riemannian manifold or locally finite polyhedral complex
(more generally, a local Lipschitz neighborhood retract).  Singular
Lipschitz $k$-chains (sometimes simply called Lipschitz $k$-chains) are formal sums
of Lipschitz maps from the standard simplex $\Delta^k$ to $X$.  As in
the cellular case, we will consider chains with integral
coefficients.  The boundary operator is defined as for singular
homology.  Rademacher's Theorem implies that a Lipschitz map is differentiable almost everywhere,
so we can define the $k$-volume of a Lipschitz map $\Delta^k\to X$ as the
integral of the Jacobian, and
we define the {\em mass} of a Lipschitz chain to be the total volume of its
summands, weighted by the coefficients.  For $k$-connected $X$, if $\CLip_{k}(X)$ 
is the set of Lipschitz $k$-chains in $X$ and $\ZLip_{k}(X)$ is the set of
Lipschitz $k$-cycles then we can define filling functions by
$$
\delta^{k;\text{Lip}}_X(a):=\mathop{\inf_{b\in \CLip_{k+1}(X)}}_{\partial b=a} \mass(b) \quad
; \qquad
\delta^{k;\text{Lip}}_{X}(\V) = \mathop{\sup_{a\in \ZLip_{k}(X)}}_{\mass a\le \V} \delta^{k;\text{Lip}}_X(a).
$$

These two definitions of Dehn functions are very similar, and if
$X$ is a polyhedral complex with bounded geometry (or if $X$ is a
space which can be approximated by such a polyhedral complex), one can show that
they grow at the same rate by using the Federer-Fleming Deformation
Theorem.  We briefly explain the notation before stating the theorem: 
we will be approximating a
Lipschitz chain $a$ by a cellular chain $P(a)$.  This may
necessitate changing the boundary, and $R(a)$ interpolates between
the old and new boundaries.  Finally, $Q(a)$ interpolates between $a$
and $P(a)+R(a)$.  Note that if $X$ is a polyhedral complex, 
then there is an inclusion $\Ccell_k(X)\hookrightarrow\CLip_k(X)$.

\begin{theorem}[Federer-Fleming \cite{FedFlem}]\label{fed-flem}\label{thm:fedflem}
  Let $X$ be a polyhedral complex with finitely many isometry types of cells.  There is
  a constant $c$ depending on $X$ such that if $a\in \CLip_k(X)$, then
  there are $P(a)\in \Ccell_k(X)$, $Q(a)\in \CLip_{k+1}(X)$, and
  $R(a)\in \CLip_k(X)$ such that   \raggedcolumns
  \begin{multicols}{2}
    \begin{enumerate}
      \item  $\|P(a)\|\le c \cdot\mass(a)$
      \item  $\|Q(a)\|\le c \cdot\mass(a)$
      \item  $ \|R(a)\|\le c \cdot\mass(\partial a)$
      \item $  \partial Q(a) = a - P(a) - R(a) $
      \item   $ \partial R(a) = \partial a - \partial P(a).$
    \end{enumerate}
  \end{multicols}
  If $\partial a \in \Ccell_k(x)$, we can take $R(a)=0$.  Furthermore,
  $P(a)$ and $Q(a)$ are supported in the smallest subcomplex of $X$
  which contains the support of $a$, and $R(a)$ is supported in the smallest
  subcomplex of $X$ which contains the support of $\partial a$.
\end{theorem}
This version of the Federer-Fleming theorem is close to the one in
\cite{wordproc}, which addresses the case that $a$ is a cycle.  

As an application, it is straightforward to prove that if $X$ is as above, then 
$\delta^{k;\text{Lip}}_{X}(\V)\asymp
\delta^{k;\text{cell}}_{X}(\V)$.   We will thus generally refer to
$\delta^{k}_{X}$ or $\delta^{k}_{G}$, using cellular or
Lipschitz methods as appropriate.  

Another (very important) application of the Federer-Fleming theorem is the isoperimetric inequality in Euclidean space~\cite{FedFlem}: if $1\le k\le d-1$, then 
$$\delta^k_{\R^d}(\V)\asymp \V^{\frac{k+1}{k}}.$$
This is extended to general
CAT(0) spaces, obtaining the same upper bound, in \cite{groFRM, Wenger}.  
We state the version we will need for our filling results; 
it describes the key properties of Wenger's construction.

\begin{prop}[CAT(0) isoperimetric inequality \cite{Wenger}] \label{prop:wenger}
  If $X$ is a CAT(0) polyhedral complex and $k\ge 1$, then the $k$-th order Dehn
  function of $X$ satisfies
  $$\delta^k_{X}(\V)\preceq \V^{\frac{k+1}{k}}.$$
  In fact, a slightly stronger condition is satisfied: there is
  a constant $m$ such that if $a\in
  \ZLip_k(X)$, then there is a chain $b\in \CLip_{k+1}(X)$ such
  that $\partial b = a$,
  $$\mass b \le m (\mass a)^{\frac{k+1}{k}},$$
  and 
  $\supp b$ is contained in a $m(\mass
  a)^{\frac{1}{k}}$-neighborhood of $\supp a$.
  \end{prop}

\subsection{Higher divergence functions}\label{sec:defineDivergence}
We will define divergence invariants $\div^k(X)$ for spaces $X$ with sufficient
connectivity at infinity.
Our goal is to study the divergence functions of groups.  We will solve 
filling problems in model spaces, such as $K(G,1)$ spaces and other cell complexes
with a geometric $G$-action.  To make this meaningful, we therefore want 
$\div^k$ to be invariant under quasi-isometries.  The somewhat complicated 
equivalence relation defined in this section 
is designed to achieve this. 

$\div^k(X)$ will basically generalize the definitions of divergence found in Gersten's 
work for $k=0$ and Brady-Farb for $k\ge 1$ \cite{gersten1, brady-farb}.  However, ours
 is not quite the same notion of equivalence.  In particular, 
ours distinguishes polynomials of different degrees, whereas
Brady-Farb identifies all polynomials into a single class by the equivalence relation used to 
define  $\div^k$.  The equivalence classes here are strictly finer than theirs.
Moreover, there is a subtle error in the definition of $\preceq$ found in Gersten's 
original paper (making the equivalence classes far larger than intended) 
that propagated into the rest of the literature.

Let $X$ be a metric space with basepoint $x_0$.  
Recall from above that a map to $X$ is $r$-avoidant if its image 
is disjoint from the ball of
radius $r$ about $x_0$.  We say that a chain (Lipschitz singular or
cellular) in $X$ is $r$-avoidant if its support is disjoint from the ball of
radius $r$ about $x_0$.
For $\rho\le 1$, we say that $X$ is $(\rho, k)$-{\em acyclic at infinity} if for 
every $r$-avoidant $n$-cycle $a$ in $X$, where $0 \le n \le k$, there exists a 
$\rho r$-avoidant $(k+1)$-chain $b$ with $\partial b =a$.  
For fixed $k$, we sometimes write $\bar \rho$ for the supremum of the values for
which $(\rho,k)$-acyclicity at infinity holds.
Note that if $X$ is $(\rho, k)$-acyclic at infinity 
for any $\rho$ then it is $k$-acyclic at infinity 
(cf.~\cite{brady-meier} for the  definition of acyclicity at infinity). The converse is false in general,
and we will discuss the special case of right-angled Artin groups in more detail 
in the next section.

For a metric space $X$, define the {\em divergence dimension} $\divdim(X)$ to be the largest
whole number $k$ such that $X$ is $(\rho,k)$-acyclic at infinity for some $0<\rho\le 1$. 
For instance, $\divdim(\R^d)=\divdim(\H^d)=d-2$.  
We will define $\div^k$ for $k\le \divdim(X)$.

The definition of $\div^k$ will be a bit special when $k=0$, so we 
deal with that case later.  Suppose first that $1\le k\le \divdim(X)$.
Given a $k$-cycle $a$, we define 
$$
\ddiv^k_\rho(a, r) := \inf \mass b,
$$
where the $\inf$ is over $\rho r$-avoidant Lipschitz $(k+1)$-chains $b$ such that $\partial b =a$.  
We then define 
$$
\ddiv^k_\rho(l, r) := \sup \ddiv^k_\rho(a, r), \qquad\qquad (k>0)
$$
where the $\sup$ is over $r$-avoidant Lipschitz $k$-cycles $a$ of mass at most $l$.

In order to see the effect of removing a ball from the space,  consider
what happens as $r$ and $l$ go to infinity simultaneously.  In the nonpositively curved
setting, the difficulty of filling spheres that arise as the intersection of a large ball 
around the basepoint with a flat of rank $k+1$ tends to be a distinguishing feature
of the asymptotic geometry (as in the results for symmetric spaces referenced above).  
These spheres have $l=O(r^k)$, and so we can obtain useful information by 
specializing to spheres whose mass is of this order.  We therefore introduce a new
parameter $\alpha$ and write
$\ddiv_{\rho,\alpha}^k(r)$ for $\ddiv_\rho^k(\alpha r^k, r)$.  
Then, formally, $\div^k(X)$ is the two-parameter 
family of functions: 
$$\div^k(X) := \{ \ddiv_{\rho,\alpha}^k (r)\}_{\alpha,\rho}, \qquad\qquad (k>0)$$
where $\alpha>0$ and $0< \rho \le \bar \rho$.

In the case $k=0$, we are filling pairs of points 
($0$-cycles) by paths ($1$-chains).  The $0$-mass of a cycle does not
restrict its diameter, so instead we require the $0$-cycle to lie on the boundary of the deleted ball.
Set 
$$\ddiv_\rho^0(r) :=  \sup_{x,y\in S_r} \inf_P |P|,$$
where the $\sup$ is over pairs of points on $S_r$ and the $\inf$ is over $\rho r$-avoidant
paths $P$ with endpoints $x$ and $y$.

In this case we get a  one-parameter family of functions of one variable:
$$\div^0(X) := \{\ddiv_\rho^0(r)\}_\rho,$$
where $0\le \rho \le \bar \rho$.

Let $F= \{ f_{\rho, \alpha}\}$ and $F'= \{ f'_{\rho, \alpha}\}$ be
two-parameter families of functions $\R^{+}\to \R^{+}$,
indexed over $\alpha>0$ and $0<\rho\le \bar\rho$.  
Then we write $F \preceq F'$ if there exist thresholds 
$0<\rho_0\le \bar \rho, \alpha_0\ge0$, and constants $L, M,A>1$ such that for all $\rho \le \rho_0$ 
and all $\alpha\ge \alpha_0$, $x>0$,
$$f_{\rho, \alpha}(x) \leq A\cdot f'_{L\rho, M\alpha}(Ax+A)+ O(x^k).$$
(Since the volume of the objects we're filling is on the order of
$r^k$, we need an $O(x^k)$ term rather than an $O(x)$ term to preserve
quasi-isometry invariance.)

\begin{figure}[ht]
\begin{tikzpicture}[scale=.8]
\foreach \a in {0,5}
{\begin{scope}
\draw (\a,0) -- (\a+3,0)--(\a+3,3)--(\a,3)--cycle;
\draw[->] (\a,0)--(\a,3.7);
\node at (\a,3.7) [left] {$\alpha$};
\node at (\a,3) [left] {$\infty$};
\draw[->] (\a,0)--(\a+3.7,0);
\node at (\a+3.7,0) [below] {$\rho$};
\node at (\a+3,0) [below] {$1$};
\end{scope}
}
\filldraw [shading=radial, inner color=gray!80, outer color=gray!30] (0,1.7)--(1.3,1.7)--(1.3,3)--(0,3)--cycle;
\draw [dashed] (1.3,1.7)--(1.3,0) node [below] {$\rho_0$};
\node at (0,1.7) [left] {$\alpha_0$};
\filldraw [shading=radial, inner color=gray!80, outer color=gray!30] (5,2.1)--(6.95,2.1)--(6.95,3)--(5,3)--cycle;
\draw [dashed] (6.95,2.1)--(6.95,0) node [below] {$L\rho_0$};
\node at (5,2.1) [left] {$M\alpha_0$};
\end{tikzpicture}
\caption{Comparison of two-parameter families of functions:  for each coordinate position
in the rectangle on the left, there is a corresponding position in the rectangle on the right.
We say $F\preceq F'$ if the functions in those positions satisfy $f\preceq f'$ over the whole
rectangle.\label{scalebox}}
\end{figure}
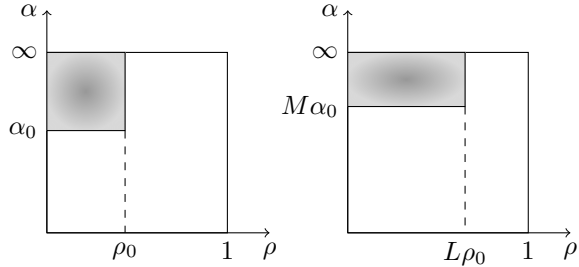

From this description it is clear that $\preceq$ is a partial order (see Figure~\ref{scalebox}).
Finally, $F \asymp F'$ if and only if $F \preceq F'$ and $F \succeq F'$.

\begin{proposition}[Quasi-isometry invariance]\label{prop:QI}
Let $X$ and $Y$ be $k$-connected cell complexes with finitely many isometry types of cells.  
If $X$ is quasi-isometric to $Y$ and $Y$ is $k$-acyclic at infinity, then $\div^k(X) \asymp \div^k(Y)$.
\end{proposition}

Quasi-isometry invariance allows us to write $\div^k(G)$ for the equivalence
class of two-parameter families $\{\div^k(X)\}$ where $X$ has a geometric $G$-action.

\begin{proof} 
If $k=0$ then this is a result of Gersten \cite{gersten1}.  Fix $k>0$ in the indicated range.  

Some technical lemmas from \cite{AlWaPr} imply the 
following:  if cell complexes $X$ and $Y$ are $k$-connected and have finitely many 
isometry types of cells, and if $X$ and $Y$ are quasi-isometric, then there are quasi-isometries 
$\varphi:X \to Y$ and $\overline\varphi:Y\to X$ that are quasi-inverses of each other 
and that are cellular and $C$-Lipschitz on the $(k+1)$-skeleton for some $C\ge 1$.
Fix such maps $\varphi,\overline\varphi$ and constant $C$.  It is furthermore possible 
to choose a constant $M$, dependent on $C$ and $k$, such that (i) the mass of any 
push-forward $\varphi_\#(\sigma)$ or $\overline\varphi_\#(\sigma)$ is at most 
$M\cdot\mass(\sigma)$ for any Lipschitz $k$- or $(k+1)$-chain in $X$ or $Y$, and 
(ii) every Lipschitz $k$-chain $a$ in $X$ is homotopic to $\overline\varphi_\#\varphi_\#(a)$ 
by a Lipschitz homotopy of mass at most $M\cdot\mass(a)$.

Let $0<\bar{\rho}\leq 1$ be the maximal value for which 
$Y$ is $(\bar{\rho},k)$-acyclic at infinity.  Let
$\rho_0=C^{-2}\bar{\rho}$.  
Let $\alpha_0=0,$ let $L=C^2,$ and 
let $M$ be as described in the previous paragraph.  Now fix $0<\rho\leq\rho_0$ and 
$\alpha\geq\alpha_0$.  We will show that $X$ is $(\rho_0,k)$-acyclic
at infinity and that
$$\ddiv^k_{\rho,\alpha}(X)\preceq\ddiv^k_{L\rho,M\alpha}(Y),$$
from which we conclude $\div^k(X)\preceq\div^k(Y)$.  A symmetric argument shows 
the other inequality, giving the desired equivalence.

Specifically, let $r>0$ be given and let $a$ be an $r$-avoidant Lipschitz $k$-cycle 
in $X$ with mass $\le \alpha r^k$.  It suffices to show that $a$ can be filled by a 
$\rho r$-avoidant Lipschitz $(k+1)$-cycle $b$ that has mass at most 
$A\ddiv^k_{L\rho,M\alpha}(Y)(r/C)$ for some constant $A>0$.

Note that the pushforward $a'=\varphi_\#(a)$ is a Lipschitz $k$-cycle in $Y$; it is 
$r/C$-avoidant and has mass at most $M\alpha r^k$.  Therefore for any 
$0<\rho'<\bar\rho(Y)$ there exists a filling $b'$ of $\varphi_\#(a)$ (that is, $b'$ is a 
Lipschitz $(k+1)$-chain) that is $(\rho' r/C)$-avoidant and that satisfies
$$\mass(b')\leq \ddiv^k_{\rho',M\alpha}(r/C).$$
Choose $\rho'=L\rho$, so that $b'$ is $C\rho r$-avoidant in $Y$.

Now consider $b''=\overline\varphi_\#(b')$.  This is a Lipschitz $(k+1)$-chain in 
$X$ that is $\rho' r/C^2=\rho r$-avoidant and that has mass at most $M\mass(b').$  
However $b''$ is not quite a filling of $a$; we know only that its boundary $a''$ is 
bounded distance from $a$.  Since $a''=\overline\varphi_\#\varphi_\#(a)$ there is 
a (Lipschitz) homotopy between $a''$ and $a$ with mass at most $M\cdot\mass(a)$; 
thus we have 
$$\ddiv^k_{\rho,\alpha}(X)(r)\leq M\ddiv^k_{L\rho,M\alpha}(Y)(r/C)+M\alpha r^k$$
and
$\ddiv^k_{\rho,\alpha}(X)\preceq \ddiv^k_{L\rho,M\alpha}(Y)$, as desired.
\end{proof}

For a function $h:\R^+\to\R^+$, we say that the family $F$ has order $h(r)$ and write 
$F\asymp h(r)$ if $F\asymp \{h(r)\}$, that is, $F$ is equivalent to 
the family that contains the same function $h(r)$ for every value of the parameters.
Then $F \preceq h(r)$ and $F \succeq h(r)$ can be defined similarly.

\begin{remark}[Remarks on comparison and equivalence] \label{compare} \ \\
\begin{enumerate}

\item Note that the statement $h(r) \preceq \div^k(X)$ is equivalent to the statement that 
there exist $\rho'_0$ and $\alpha'_0$ such that $h(r) \preceq \ddiv_{\rho, \alpha}(r)$
for all $\rho \le \rho'_0$ and $\alpha \ge \alpha'_0$.   (Here $\rho_0'=\rho_0 L$ and $\alpha_0' = \alpha_0 M$ in the definition of $\preceq$.)

Similarly, note that if  $\rho$ or $\alpha$ is decreased, then the value of 
$\ddiv^k_{\rho, \alpha}(r)$ decreases.  Thus in order 
to establish that $\div^k(X) \preceq h(r)$, it suffices to show that there exist $\rho_0 \le 1$ 
and $\alpha_0 \ge 0$ such that $ \ddiv_{\rho_0, \alpha}(r) \preceq h(r)$ for all 
$\alpha \ge \alpha_0$.

Taken together, these give sufficient criteria to establish that $\div^k(X)\asymp h(r)$, 
as shown in Figure~\ref{block-compare}.
However, for
 a particular $X$ there is no guarantee that $\div^k(X)\asymp h(r)$ for any $h$.

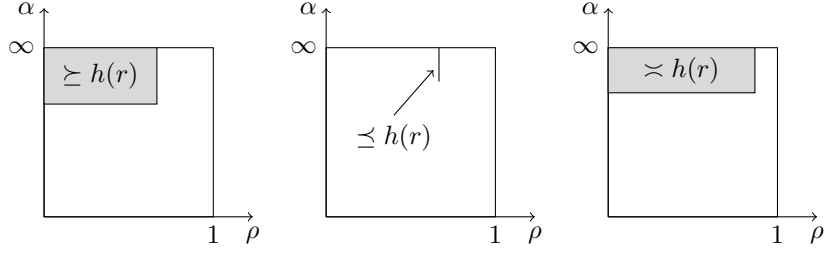
\begin{figure}
\begin{tikzpicture}[scale=3/4]

\foreach \a in {0,5,10}
{\draw (\a,0) -- (\a+3,0)--(\a+3,3)--(\a,3)--cycle;
\draw[->] (\a,0)--(\a,3.7);
\node at (\a,3.7) [left] {$\alpha$};
\node at (\a,3) [left] {$\infty$};
\draw[->] (\a,0)--(\a+3.7,0);
\node at (\a+3.7,0) [below] {$\rho$};
\node at (\a+3,0) [below] {$1$};
}

\draw [fill=gray!30] (0,3) rectangle (2,2);
\node at (1,2.5) {$\succeq h(r)$};
\draw (7,3) -- (7,2.4);
\draw [->](6.2,1.8) node [below] {$\preceq h(r)$}--(6.9,2.6);
\draw [fill=gray!30] (10,3) rectangle (12.6,2.2);
\node at (11.3,2.6) {$\asymp h(r)$};

\end{tikzpicture}
\caption{The box diagrams show sufficient criteria to check that $\div^k(X)$ compares to the 
function $h(r)$ by $\succeq$, $\preceq$, and $\asymp$,
respectively, as described in Remark~\ref{compare}(1).
\label{block-compare}}
\end{figure}

\item Every family $F$ satisfies $F\succeq r$, since all sublinear functions are equivalent
to the linear function $r$.  In particular $\div^0(G)\succeq r$ for any group $G$.
The following (true) statement is slightly stronger in two ways: 
if $G$ is a finitely generated infinite group, 
then for every $\rho$ we have $\ddiv^0_{\rho}(r) \ge 2r$ for all $r>0$.

\item For CAT(0) spaces with extendable geodesics, one 
sees only one $\asymp$ class of functions 
in $\div^0(X)$.  On the other hand, for $k\ge 1$ the family
$\div^k(X)$ often contains functions from multiple $\asymp$ classes, as we will see in 
the example of orthoplex groups in \S\ref{sec:BB}.

\end{enumerate}
\end{remark}

For groups that are direct products, it is easy to see that $\div^0$ is exactly linear;
we will show this below (Lemma~\ref{lem:LinearUpper}).

\subsection{Example:  Euclidean space} \label{sec:euclidean example}

Constructing avoidant fillings is sometimes difficult; removing a ball
of radius $r$ from a space breaks its symmetry, making it harder to
apply methods from group theory.  One method of constructing avoidant
fillings is to first construct a filling in the entire space, then
modify that filling to be avoidant.  In this paper, we modify fillings
using maps $X\to X\setminus B_r$; we call these pushing maps.

Our constructions generally follow the following outline:
given an avoidant $k$-cycle $a$ in $X$, we will find a filling $b$
(typically not avoidant) of $a$ and a pushing map $X\to X\setminus
B_r$, where $B_r$ is the ball of radius $r$.  This map generally has
singularities, and we use techniques from geometric measure theory to
move $b$ off these singularities.

The basic example to consider is $\R^d$, where one has the ``pushing'' map
given by radial projection to $\R^d\setminus B_r$, namely
\begin{equation}\label{eq:euclidean-pushing}
\pi_r(v)=\begin{cases}
  r \frac{v}{\|v\|}, &  v\in B_r \\
  v, &  v\notin B_r.
\end{cases}
\end{equation}
This map is undefined at $0$, but if the filling has dimension $<d$, it can be 
perturbed to miss the origin, and the Federer-Fleming Deformation Theorem 
(Theorem~\ref{thm:fedflem}) can be used to control the volume.  The theorem
allows us to approximate singular $k$-chains in $X$ by cellular
$k$-chains in $X$, and if the $k$-skeleton of $X$ misses the
singularity, then so will the approximation.

We will prove that filling an avoidant cycle by an avoidant
chain is roughly as hard as filling a cycle by a chain.  Specifically,
we will show the following.
\begin{prop}[Euclidean bounds]\label{prop:eucBounds}
  Let $d$ be a positive integer and let $1\le k\le \divdim(\R^d)=d-2$.  
  There is a constant $c$ depending only on the dimension $d$ 
  such that if $r,\V \ge 0$ and $a$ is an $r$-avoidant $k$-cycle in
  $\R^d$ of mass at most $\V$, then there is an $r$-avoidant $(k+1)$-chain $b$ such that
  $\partial b=a$ and
  $$\mass b \le c\V^{\frac{k+1}{k}}.$$
  Further, there is a constant $c'$ depending only on $d$ and
  an $r$-avoidant $k$-cycle $a$ with mass $\V$
  such that for every chain $b$
  with $\partial b=a$, 
  $$\mass b \ge  c' \V^{\frac{k+1}k}.$$
 \end{prop}
This implies, in particular, that for $0\le k \le d-2$,
$\div^k(\R^d)\asymp r^{k+1}.$
As $\R^d$ is a model space for $\Z^d$, Proposition~\ref{prop:QI} gives
$\div^k(\Z^d) \asymp r^{k+1}.$

\begin{proof}
  Wenger's work (Prop.~\ref{prop:wenger}, though cf.\ Federer and
  Fleming \cite{FedFlem} in the Euclidean case) implies that there exists
  an $m>0$ independent of $a$ and there exists a chain $b$ such that $\partial b=a$,
  $\mass b\le m \V^{(k+1)/k}$, and $\supp b$ is contained in a
  $m\V^{1/k}$-neighborhood of $\supp a$.  We will modify this to find
  an avoidant chain, using different arguments when $\V\preceq r^k$
  and when $\V\succeq r^k$.  

  First, set $c_0=(2m)^{-k}$ and note that if $\V\le c_0r^k$ then $b$ is 
  $r/2$-avoidant.  In this case $b'=(\pi_r)_\sharp(b)$ fills $a$, is $r$-avoidant, 
  and satisfies
  $$\mass b'\le 2^k\mass b\le 2^km \V^{(k+1)/k},$$
  so the conclusion of the Proposition holds.

  We can thus assume that $\V\ge c_0 r^k$.  We will show the
  proposition when $r=1$, and then use scaling to prove the general
  case.  Let $a$ be a $1$-avoidant Lipschitz $k$-cycle of mass
  $\V\ge c_0$, and $b$ be its filling as above.  We will
  approximate $b$ by a cellular chain, then ``push'' it out of the
  $1$-sphere.

  Let $\tau$ be a grid of cubes of side length $\frac{1}{2 d}$
  translated so that the center of one of the cubes lies at the
  origin.  
  Let $P(b)$, $Q(b)$, and $R(b)$ be as in
  Federer-Fleming, so that $P(b)$ is a chain in $\tau^{(k+1)}$ approximating $b$, and 
  $$
  \partial R(b)=\partial b-\partial P(b)=a-\partial P(b).
  $$
  Each cell of $\tau$ has diameter at most $1/2$, so the
  smallest subcomplex of $\tau$ containing the support of $a$ is
  $1/2$-avoidant.  It follows that $R(b)$ is $1/2$-avoidant. 
  Since any $k$-cell of $\tau$ is $1/4d$-avoidant, so is
   $P(b)$.
   Thus
  $b':=R(b)+P(b)$ is a $1/4d$-avoidant filling of $a$.  Further, there is a
  constant $c_1$, which comes from Federer-Fleming and depends only on $d$, 
   such that
   $$\mass b'\le c_1(\mass a + \mass b)\le
   c_1(\V+m\V^{\frac{k+1}{k}})\le c_2\V^{\frac{k+1}{k}}$$ 
   for some $c_2$, where the last bound uses the lower bound on $\V$.

  Pushing $b'$ forward under the radial pushing map $\pi_1$ from \eqref{eq:euclidean-pushing}, 
  we get a chain $b'':=(\pi_1)_\sharp(b')$.  This is a $1$-avoidant filling of $a$.
  Furthermore, since $b'$ is $1/4d$-avoidant, $\pi_1$ is $4d$-Lipschitz
  on the support of $b'$, so there is a constant $c$ such that
  $$\mass b''\le (4d)^{k+1}\mass b' \le c \mass a^{\frac{k+1}{k}},$$
  as desired.

  Now, return to the case of general $r$.  
  Let $s_t:\R^d\to \R^d$ be the
  homothety $v\mapsto tv$.  If $\gamma$ is a Lipschitz $k$-chain, then
  $\mass s_t(\gamma)=t^k\mass \gamma.$

  If $a$ is an $r$-avoidant Lipschitz $k$-cycle of mass $\V$, then $a_1={s_{r^{-1}}}_\sharp(a)$
  is $1$-avoidant, and $\mass a_1=r^{-k}\V\ge c_0$ (since $\V\le c_0 r^k$ is done already).
  By the argument above,
  there is a $1$-avoidant $(k+1)$-chain $b_1$ filling $a_1$ such that
  $$\mass b_1 \le c(\mass a_1)^{\frac{k+1}{k}} = c r^{-(k+1)}\V^{\frac{k+1}{k}}.$$
  Rescaling this by letting $b=(s_{r})_\sharp(b_1)$, we find that $b$
  is an $r$-avoidant filling of $a$ and
 $\mass b \le c\V^{\frac{k+1}{k}}.$ 
 
 For the second statement, simply take $a$ to be a round sphere far from $B_r$.  
 Then the estimate is just the isoperimetric theorem for $\R^d$.
 \end{proof}

Thus the best avoidant fillings have roughly the same volume as the most efficient 
(not necessarily avoidant) fillings.  

Much of the rest of this paper will be dedicated to generalizing this
technique to right-angled Artin groups.  These groups act on a complex
$X$ which consists of a union of flats, and as with $\R^d$, we will
construct avoidant fillings by using a pushing map to modify non-avoidant fillings.  
The pushing map is singular in the sense that it cannot
be defined continuously on all of $X$, but as with $\R^d$, we will
delete small neighborhoods of the singularities, enabling us to define
the pushing map continuously on a subset of $X$.  Because $X$ generally
has more complicated topology than $\R^d$, the pushing map has more
singularities, and these singularities lead to larger bounds on
$\div^k$.


\section{Background on right-angled Artin groups}

In this section we will introduce some of the key background on RAAGs.
We refer the reader to ~\cite{charney} for a more complete treatment.

A right-angled Artin group is a finitely generated group given by a
presentation in which all the relators are commutators of generators.
A RAAG can be described by a graph which keeps track of which pairs of
generators commute, and the structure of this graph affects the
geometry of the group and its subgroups.  Let $\Gamma$ be a finite graph
with no loops or multiple edges, and with vertices labeled $a_1,\ldots, a_n$.  
The right-angled Artin group based on $\Gamma$ is the group
$$A_\Gamma:=\langle a_1,\ldots, a_n \mid R\rangle, \qquad \hbox{\rm with relators} \quad
R=\{[a_i,a_j]\mid \hbox{\rm $a_i$, $a_j$  connected by an edge of $\Gamma$}\}.$$
We call $\Gamma$ the
{\em defining graph} of $\Agamma$.  Let
$L$ be the flag complex of $\Gamma$; that is, the simplicial complex with the same vertex 
set as $\Gamma$, and in which a set $S$ of vertices spans
a simplex if and only if every pair of vertices of $S$ is connected by
an edge of $\Gamma$.  

The group $A_\Gamma$ acts freely on a CAT(0) cube-complex $X_\Gamma$,
defined as follows.  Let $Y$
be a subcomplex of the torus $(S^1)^n$, where the circle $S^1$ is given a cell structure 
with one $0$-cell and one $1$-cell, and each $S^1$ factor corresponds to a vertex 
$a_i$.  Thus $Y$ is a cube complex with one vertex.  A $d$-cell $\sigma$ of $(S^1)^n$ 
is contained in $Y$ if and only if the vertices corresponding to the $S^1$ factors of 
$\sigma$ span a simplex in $L$.  Then $X_\Gamma$ is the universal cover of $Y$.

Since $Y$ has one vertex, all the vertices of $X_\Gamma$ are in the
same $\Agamma$-orbit.  We pick one of the vertices of $X_\Gamma$ as a basepoint,
which we denote $e$, and identify the element $a\in \Agamma$ with the
vertex $a\cdot e$ of $X_\Gamma$.  We will often refer to elements of
$\Agamma$ and vertices of $X_\Gamma$ interchangeably.  Similarly, the
edges of $Y$ correspond to the generators of $A_\Gamma$, and we say
that an edge of $X_\Gamma$ is {\em labeled} by its corresponding
generator.

Since each cell in $Y$ is part of a torus, each cell of $X_\Gamma$ is
part of a flat.  Typically, each cell is part of infinitely many
flats, but we can use the group structure to pick a canonical one.  If
$\sigma$ is a $d$-cube in $X$, its edges are labeled by $d$ different
generators; if $S$ is this set of labels, and if $v\in A_\Gamma$ is a
vertex of $\sigma$, then the elements of $S$ commute and generate an undistorted
copy of $\Z^d$, which we denote $A_S$.  This subgroup spans a
$d$-dimensional flat through the origin, and the translate $v\cdot A_S$ 
spans a flat containing $\sigma$, which we call the {\em standard flat} 
containing $\sigma$ and denote by $F_\sigma$.

The link of a vertex of $X_\Gamma$, which we denote $S(L)$, is the
union of the links of all the standard flats.  It has two vertices for
each generator $v$ of $A_\Gamma$, one corresponding to $v$ and one to
$v^{-1}$.  We will denote the vertex in the $v$ direction by
$\hat{v}=+\hat{v}$ and the one in the $v^{-1}$ direction by
$-\hat{v}$.  Furthermore, if $v_1,\dots,v_d$ are the vertices of a
simplex $\Delta$ of $L$, that simplex corresponds to a $d$-dimensional
standard flat containing $x$.  The link of this flat is an \emph{orthoplex} (i.e.,
the boundary of a cross-polytope; in the case that $d=3$, it is an octahedron), 
so $\Delta$ corresponds to
$2^d$ simplices in $S(L)$, each with vertices $\pm \hat{v}_1,\dots, \pm
\hat{v}_d$.  If $L$ has $m$ vertices, labeled $v_1,\dots, v_m$, then
$S(L)$ contains $2^m$ copies of $L$ with vertices $\pm v_1,\dots, \pm
v_m$; we call these {\em signed copies} of $L$.

We will use two metrics on $X_\Gamma$.  The first metric on
$X_\Gamma$, with respect to which it is CAT(0), is the Euclidean (or
$\ell^2$) metric on each cube, extended as a length metric to
$X_\Gamma$.  (That is, the distance between two points is the infimal
length of a path connecting them, where length is measured piecewise
within each cube.)  

The second metric restricts instead to the $\ell^1$ metric on each cube, and is 
extended as a length metric from the cubes to the whole space.  This has the
property that its restriction to the one-skeleton of $X_\Gamma$
is the word metric on a Cayley graph for $A_\Gamma$.  
We will mainly use the
$\ell^1$ metric to define balls and spheres in $X_\Gamma$ which
coincide with balls and spheres in $A_\Gamma$.  The notation $B_r(x):= \{ y\in X :
d_{\ell^1}(x,y) < r\}$ will denote the open $\ell^1$ ball and
$\bar{B}_r(x)$, $S_r(x)$ will denote the closed ball and sphere,
respectively, so that $B_r(x) \sqcup S_r(x) = \bar{B}_r(x)$.  When
there is no center specified for a ball, it is taken to be centered at
the basepont $e$.  Note that all words whose reduced spellings have length $r$ are
vertices in $S_r$.

 As an illustration, $\Z^3$ is a RAAG, and the corresponding $X_\Gamma$ is $\R^3$, with the 
 structure of a cube complex. 
The sphere of radius $r$ in the $\ell^1$ metric is a Euclidean octahedron with equilateral triangle 
faces.  All vertices corresponding to group elements of length $r$ 
in the word metric lie on this sphere.

Recall that RAAGs themselves, being CAT(0) groups, have at worst Euclidean Dehn 
functions (Proposition~\ref{prop:wenger}).  To find bigger Dehn functions, one must look 
at subgroups such as those defined in the next section.

We will study divergence functions for RAAGs below, so we remark that the 
divergence dimension can be read off of the defining graph.
Brady and Meier showed that 
the group $\Agamma$ is $k$-acyclic at infinity if and only if the link $S(L)$ is $k$-acyclic.
In fact, their construction shows that $k$-acyclicity of the link is equivalent to 
$(1,k)$-acyclicity at infinity
of the group (and therefore $(\rho,k)$-acyclicity at infinity for all $0<\rho\le 1$).
Thus,
$\divdim(\Agamma)$ is the largest $k$ such that $S(L)$ is $k$-acyclic.

\subsection{Bestvina-Brady groups}

Let $h:A_\Gamma\to \Z$ be the homomorphism which sends each generator
to $1$; we call $h$ the {\em height function} of $A_\Gamma$.  Let
$H_\Gamma:=\ker h$; a group $H_{\Gamma}$ constructed in this fashion is called a
Bestvina-Brady group.  These subgroups were studied by Bestvina and
Brady in \cite{bestvina-brady}, and provide a fertile source of
examples of groups satisfying some finiteness properties but not
others.  Brady \cite{bradybook} showed that there are graphs $\Gamma$
such that $H_\Gamma$ has a quartic ($\V^4$) Dehn function, and Dison
\cite{dison} recently showed that this is the largest Dehn function
possible, that is, the Dehn function of any Bestvina-Brady group is at
most $\V^4$.  We will generalize Dison's result to higher-order Dehn
functions in Section~\ref{sec:BB} below.

Abusing notation slightly, 
let $h:X_\Gamma\to \R$ also denote the usual height map defined by linear extension of the homomorphism above;
it is a Morse function on $X_\Gamma$, in the sense of \cite{bestvina-brady}.
Let $Z_\Gamma=h^{-1}(0)$ be the zero level set.  The action of
$A_\Gamma$ on $X_\Gamma$ restricts to a geometric action of $H_\Gamma$ on
$Z_\Gamma$.  The topology of $Z_\Gamma$ is closely related to
$\Gamma$; indeed, if $L$ is the flag complex corresponding to
$\Gamma$, then $Z_\Gamma$ contains infinitely many scaled copies of
$L$ and is homotopy equivalent to a wedge of infinitely
many copies of $L$ \cite{bestvina-brady}.
  In particular, if $L$ is
$k$-connected, then $Z_\Gamma$ is also $k$-connected, so $H_\Gamma$ is
type $F_{k+1}$.

\subsection{Tools for RAAGs}\label{RAAGtools}
We introduce several basic tools:  the {\em orthant} associated to a cube in the complex
$X=\Xgamma$, a related {\em scaling map} on $X$, and an {\em absolute value map} on $X$.

Fix attention on a particular $d$-cube $\sigma$ in $X$ 
and let $v$ be its closest vertex to the origin.  The vertices of the standard flat $F_\sigma$ 
correspond to a coset $vA_S$ where $S$ is the set of labels on edges of $\sigma$.
For each $a_i\in S$ let $\gamma_i$ be the geodesic ray in $F_\sigma$ that starts at $v$, 
traverses the edge of $\sigma$ labeled $a_i$ in time one, and continues at this speed
inside $F_\sigma$, so that
$\gamma_i(n)=va_i^{\pm n}$, with the sign fixed once and for all depending on whether 
$v$ or $va_i$ is closer to the origin.
We take $\Orth_\sigma$ to be the flat orthant within $vA_S$
spanned by the rays $\gamma_i$, 
so that the cube $\sigma$ itself is contained in $\Orth_\sigma$, and $v$ is its extreme point.
Note that if $\tau$ is
a face of $\sigma$, then $\Orth_{\tau}\subset \Orth_\sigma$ as an orthant with appropriate
codimension.  In particular, if $\tau$ is an edge, then $\Orth_\tau$ is a ray 
starting at one endpoint of $\tau$ and pointing ``away'' from $e$.

Next we define a {\em scaling map} $s_r: S(L)\to \Xgamma$.
The sphere $S_1$ (the unit sphere in the $\ell^1$ metric) is homeomorphic to $S(L)$; we
associate points of $S(L)$ with points of $S_1$.  Because of the
abundance of flats in $X_\Gamma$, these correspond to canonically extendable directions in
$X_\Gamma$, as follows.  If $x\in S_1$, then $x$ is in some maximal cube $\sigma$
corresponding to commuting generators;
we define $\gamma_x:[0,\infty)\to X_\Gamma$ to be the unique geodesic ray
in $F_\sigma$ that is based at the identity, goes through $x$, and is parametrized by 
arc length in the $\ell^1$ metric.  For instance, if $x$ is a vertex corresponding to a generator $a$, then 
$\gamma_x$ is a standard ray along edges labeled $a$, so that $\gamma_x(n)=a^n$
for $n=0,1,2,\ldots$.  Note that the map $x\mapsto \gamma_x$ is continuous.  The 
scaling map is defined by $s_r(x)=\gamma_x(r)$.

Finally we define the {\em absolute value} map.  Given an element $g\in
A_\Gamma$, let
$w=a_{i_1}^{\pm 1}\dots a_{i_r}^{\pm 1}$
be a geodesic word representing $g$.  Then the absolute value of $g$ is given by 
$$\abs(g)=a_{i_1}\dots a_{i_r}.$$
We claim this is well-defined; indeed, if $w$ and $w'$ are two
geodesic words representing $g$, then $w$ can be transformed to $w'$
by a process of switching adjacent commuting letters \cite{herm-meier}.  If
$a_{i_j}^{\pm 1}$ and $a_{i_{j+1}}^{\pm 1}$ commute, then so do
$a_{i_j}$ and $a_{i_{j+1}}$, so the choice of geodesic spelling does not
affect $\abs(g)$.

The absolute value map preserves adjacencies.  If $g_1$
and $g_2$ are adjacent in the Cayley graph of $A_\Gamma$, then since
$A_\Gamma$ has no relations of odd length, we may assume that
$|g_1|+1=|g_2|$.  Let $a$ be a generator such
that $g_2=g_1a^{\pm 1}$.  If $w$ is a geodesic word representing $g_1$, then
$w a^{\pm 1}$ is a geodesic word representing $g_2$, so
$\abs(g_2)=\abs(g_1)a$.

Like the height function $h$, $\abs$ can be extended to $X_\Gamma$ by extending linearly
over each face.  This extension is $1$-Lipschitz and satisfies the
property that $h(\abs(x))=|x|=|\abs(x)|$ for all $x\in
X_\Gamma$; in other words, $\abs$ maps the $r$-sphere $S_r$ into the
coset $h^{-1}(r)$ of $H_\Gamma$.  Furthermore, $\abs$ is idempotent; if $h(x)=|x|$,
then $\abs(x)=x$.

\section{Pushing maps}\label{sec:pushing}

Throughout this section, let $\Gamma$ be the defining graph of a RAAG,
let $X=X_\Gamma$, $A=A_\Gamma$, and $Z=Z_\Gamma$.

In \S\ref{sec:euclidean example}, we constructed avoidant fillings
in $\R^d$ by using a pushing map $\R^d \setminus\{0\}\to \R^d\setminus B_r$; in this
section, we will construct pushing maps for general RAAGs.  Here the pushing maps 
become more complicated, because typically there is branching at the vertices.  
This has two implications:  first, there are many more singularities to work around,
and second, one needs to be careful in computing the amount by which the map 
distorts volumes.  In contrast to the situation in $\R^d$, where
avoidant fillings are roughly the same size as ordinary fillings,
avoidant fillings in RAAGs (when they exist) may be much larger.

We will define two pushing maps:  heuristically, one pushes radially to
$X\setminus B_r$ from the basepoint, and the other pushes along the height
gradient to the $0$-level set $Z$.  The constructions are similar; in both cases we delete
neighborhoods of certain vertices, put a cell structure on the remaining space, and then 
define a map cell by cell to the target space.  The Lipschitz constants of the maps
produce upper bounds on the volume expansion of fillings.

Given $X=\Xgamma$, we define modified spaces $X_r$ and $Y$, which have
$\ell^1$-neighborhoods of some vertices removed.  Let
$$X_r = X\setminus \bigcup_{v\in B_r} B_{1/4}(v) \qquad \hbox{\rm and} \qquad  Y
= X \setminus \bigcup_{v\notin Z} B_{1/4}(v).$$ 
Endow $X_r$ and $Y$ with their length metrics from the $\ell^1$ metric on $X$.
In a moment we will describe cell structures on the spaces $X_r$ and $Y$.
Recall that $B_t$ denotes the open ball, so that a modified space 
$X\setminus B_t(v)$ includes the boundary $S_t(v)$.

\begin{theorem}[Radial pushing map]\label{pushingmap}
For $r>0$ there are  Lipschitz retractions
$$\P_r: X_r \to X \setminus B_r$$  
which are $O(r)$-Lipschitz.  That is, there is a constant $c=c(\Gamma)$
such that $\P_r$ is $(cr+c)$-Lipschitz for each $r$.
\end{theorem}

\begin{proof}
We first give $X_r$ the structure of a polyhedral complex.  If $\sigma$ is a cell of $X$ 
that intersects $X_r$, then $\sigma \cap X_r$ is a cube or a truncated cube, and we
let $\sigma \cap X_r$ be a cell of $X_r$.  We call such a cell an \emph{original} face
of $X_r$.  The boundary of an original face consists of other original faces
(cubes and truncated cubes) as well as some simplices arising from the
truncation; these simplices are also cells of $X_r$ and we call them \emph{link}
faces.  Note that every link face is a simplex in a translate of $S_{1/4}$, which
we identify with $S(L)$.
If $\tau$ is an original face of $X_r$, so that $\tau=\sigma\cap X_r$ for some cube 
$\sigma$ of $X$, then we define $\Orth_{\tau}=\Orth_{\sigma}$.

\begin{figure}[ht]
\includegraphics[width=3in]{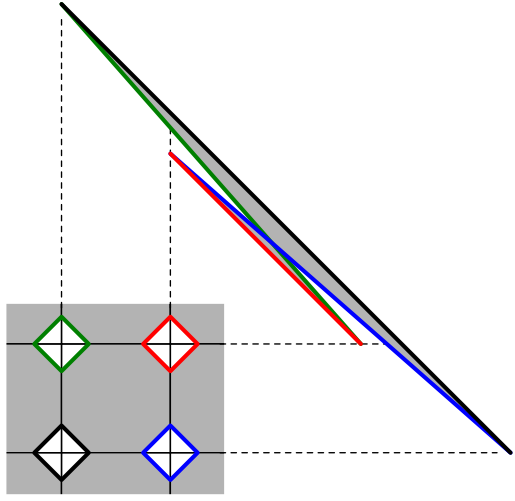}
\caption{This figure shows a portion of $X$ with $X_r$ (a union of
  truncated squares) shaded.  The vertical and horizontal edges are {\em original edges} 
and the thick diagonal lines are {\em link edges}.
The map $\P_r$ sends each original edge to a point on $S_r$, and sends 
each link edge to a line segment in $S_r$ (vertices perturbed to avoid
overlaps).  The boundary of each octagonal cell is sent to a loop of
length $\asymp r$ with zero area.
}
\label{stopsign}
\end{figure}

The pushing map $\P_r$ is the identity on $X \setminus B_r$.  
To define it on $X_r \cap B_r$, we will first define a map on the
original edges, then extend it linearly to the link faces, and finally extend
inductively to the remaining original faces.  We will require that
$\P_r(\tau\cap B_r)\subset \Orth_\tau \cap S_r$ for every original
face $\tau$ of $X_r$.  Note that $\Orth_\tau \cap S_r$ is the intersection of
$\Orth_\tau$ with a hyperplane.

We first consider the original edges of $X_r$.  If $\tau$ is an
original edge, it is part of a ray $\Orth_\tau$ (in $X$)
traveling away from the origin; we push all
  its points along $\Orth_\tau$ until they hit $S_r$, setting
  $\P_r(\tau) = \Orth_\tau\cap S_r$ (so that the image of 
  each such edge is a single point on $S_r$).
  To be explicit, suppose $\tau$ comes from an edge with endpoints
  $v$ and $va$, with $|v|,|va|\le r$.  Then for a point $x\in\tau$, 
$$ \P_r(x)= \begin{cases} va^{r-|v|} & \text{if $|va|>|v|$}; \\ 
    va^{|v|-r}  & \text{if $|va|<|v|$}.
  \end{cases}
  $$

Then the image of each original edge is a point in $S_r$ and the
images of adjacent edges are separated by distance $\preceq r$.  If
two edges of $X$ are adjacent, points on the corresponding original
edges in $X_r$ are separated by distance at least $1/4$, so the map is
Lipschitz on the edges with constant $\preceq r$.

If $\sigma$ is a link face, and $\tau$ is an original face which contains
$\sigma$, then $\P_r$ sends the vertices of $\sigma$ to points in
$\Orth_\tau \cap S_r$.  We can extend $\P_r$ linearly to
$\sigma$: every point $x$ of $\sigma$ is a unique convex combination
of its vertices, so we define the image $\P_r(x)$ to be the same
convex combination of the images of the vertices.  This is clearly
independent of our choice of $\tau$, so the map is well-defined.
Since incident edges have images no more than $O(r)$ apart, this
extension is also $O(r)$-Lipschitz.

It remains to define $\P_r$ on the rest of the original faces.  We
will proceed inductively on the dimension of the faces, using edges as
the base case.  
Recall that if $K\subset \mathbb{R}^d$ is a convex subset of Euclidean space and
  $f:S^m\to K$ is a Lipschitz map, we can extend $f$ to the ball $D^{m+1}$
via
  $$g(t,\theta)=f(x_0)+t (f(\theta)-f(x_0)),$$
  where $x_0$ is a basepoint on $S^m$ and $t\in [0,1]$,
  $\theta\in S^m$ are polar coordinates for $D^{m+1}$.  This extension
  has Lipschitz constant bounded by a multiple of $\Lip(f)$.

Each $d$-dimensional original face $\tau$ of $X_r$ is of the form
$\tau=\sigma\cap X_r$ for some $d$-cube $\sigma$ of $X$.  We may assume by
induction that $\P_r$ is already defined on the boundary of $\tau$.  
Define $\P_r$ to be the identity on $\tau\setminus B_r$.  The remaining 
part, $\tau'=\tau\cap B_r$,  is isometric to one of finitely many polyhedra, so
it is bilipschitz equivalent to a ball $D^d$, with uniformly bounded Lipschitz constant.
  
Since the boundary of $\tau'$ is mapped to a convex subset of a flat, namely
$\Orth_\sigma \cap S_r$, and the Lipschitz constant of $\P_r \big|_{\partial \tau'}$
is $O(r)$, we can extend $\P_r$ to a Lipschitz map sending $\tau'$ to
$\Orth_\sigma \cap S_r$ which is again $O(r)$-Lipschitz, with the constant enlarged 
by a factor depending on the dimension.
\end{proof}

A similar pushing map on RAAGs can be used to find bounds on higher-order
fillings in Bestvina-Brady groups.  We proceed slightly differently
here; instead of defining a map on the original edges and extending it
to the rest of the space, we start with a map on the link faces and
extend it.

Recall that the vertices of $S(L)$ are labeled $\pm \hat{a}_1,\dots,
\pm \hat{a}_d$, where the $a_i$ are generators of $A$.  Let
$\pos:S(L)\to S(L)$ be the simplicial map satisfying $\pos(\pm
\hat{a}_i)=\hat{a}_i$.  The image of this map is the copy of $L$
inside $S(L)$ whose vertices all have positive signs; we denote this
{\em positive link} by $S(L)^+$.  Similarly define $\neg: S(L)\to
S(L)$ so that $\neg(\pm \hat{a}_i)=-\hat{a}_i$ and define the {\em negative
link} $S(L)^-$ to be its image.

\begin{thm}[Height-pushing map]\label{pushing2}
  There is an $H$-equivariant retraction
  $$\QQ: Y \to Z,$$
  such that the Lipschitz constant of $\QQ$ grows linearly with distance
  from $Z$.  That is, there is a uniform constant $c=c(\Gamma)$ such that the 
  restriction of $\QQ$ to $h^{-1}([-t,t])$ is $(ct+c)$-Lipschitz.
\end{thm}

\begin{proof}
Here, instead of mapping an original face $\tau$ into
$\Orth_\tau \cap S_r$, we map it to $F_\tau \cap Z$,
where $F_\tau$ is the standard flat containing $\tau$.

Like $X_r$, the space $Y$ inherits the structure of a polyhedral complex from
$X$.  The cells of $Y$ are either original faces or link faces; the union of the link 
faces is a union of translates of $S_{1/4}$, which we identify with $S(L)$.
Thus a link vertex is denoted by $v\cdot(\pm \hat{a})$
for some $v\in A$ and some generator $a$ of $A$.

We construct a map on the link faces, then extend.  Each copy of $S(L)$ in the boundary of $Y$ is essentially the set of
directions at some vertex of $X$.  We construct a map $v \cdot S(L)\to Z$ by
flipping each direction, if necessary, to point towards $Z$, and then pushing along 
standard rays in $X$.  That is, if $v\in
A\setminus H$ and
$x\in S(L)$ we define
$$\QQ(v\cdot x)=\begin{cases} 
  v\cdot s_{|h(v)|}(\neg(x)) & \text{ if $h(v)>0$} \\
  v\cdot s_{|h(v)|}(\pos(x)) & \text{ if $h(v)<0$}.
\end{cases}$$
It is easy to check that the image of this map lies in $Z$.
Furthermore, the Lipschitz constant of this map restricted to $v\cdot
S(L)$ is $4|h(v)|$, and if $\sigma$ is a link face of $Y$ contained in
an original face $\tau$, then $\QQ(\sigma)\subset F_\tau \cap Z$.

We have defined the map on link faces of $Y$, and we extend to the
rest of $Y$ by the same inductive procedure as in
Theorem~\ref{pushingmap}, obtaining an Lipschitz constant of order $t$
on $h^{-1}([-t,t])$.
\end{proof}

\begin{remark}[Signed copies of $L$]\label{rem:signedL}
A fact that will be useful in the sequel is that if $L'$ is a signed
copy of $L$ in $S(L)$ (that is, an isomorphic copy of $L$ with some of the
vertices of $L$ replaced by their negatives), then $\QQ(L')$ is a
scaled copy of either $\pos(L)$ or $\neg(L)$; recall that Bestvina and
Brady proved that $Z$ is a union of infinitely many scaled copies of
$L$ \cite{bestvina-brady}.
\end{remark}

In \S\ref{sec:BB}--\ref{sec:higherdiv}, we will start with arbitrary $(k+1)$-dimensional fillings, and 
approximate them by fillings in the $(k+1)$-skeleton.  In order to use the pushing maps to construct fillings that are either $r$-avoidant or in $Z$, we need to extend the maps $\P_r$ and $\QQ$ to the $(k+1)$-skeleta of the deleted balls.   The connectivity hypotheses in the next two lemmas will be satisfied in the applications. 

\begin{lemma}[Extended radial pushing map]\label{lem:P}
If $S(L)$ is $k$-connected, the radial pushing map $\mathcal P_r$ can 
 be extended to an $O(r)$-Lipschitz map
  $$\P_r:X_r\cup X^{(k+1)}\to X \setminus B_r.$$
\end{lemma}
\begin{proof}
 This requires extending 
  $\P_r$ to the $(k+1)$-skeleta of the removed balls.
  Note that $B_{1/4}$ is the cone over $S(L)$, so it is a simplicial
  complex in a natural way.  We will define a retraction
  $$\rho: B_{1/4}^{(k+1)}\cup S(L) \to S(L)$$
  and use it to extend $\P_r$.  We construct $\rho$ by extending the map
  $\id_{S(L)}$ to $B_{1/4}^{(k+1)}$.  Since $S(L)$ is $k$-connected,
  there is no obstruction to constructing such an extension, and since
  $S(L)$ is a finite simplicial complex (and thus a compact Lipschitz
  neighborhood retract) one can choose it to be Lipschitz.

  Then if $v$ is a vertex in $B_r$, we can extend $\P_r$ to $v\cdot
  B_{1/4}^{(k+1)}$ by letting $\P_r(v\cdot x)=\P_r(v\cdot \rho(x))$.  This
  extension may increase the Lipschitz constant, but
  we still have $\Lip(\P_r)=O(r)$.
\end{proof}

\begin{lemma}[Extended height-pushing map]\label{lem:Q}
 If $L$ is $k$-connected, the height-pushing map $\QQ$ can be extended to an $O(r)$-Lipschitz map
  $$\QQ:Y \cup X^{(k+1)}\to Z.$$
  \end{lemma}
  
  \begin{proof}
  As before, it suffices to extend $\QQ$
  over the $(k+1)$-skeleta of the removed balls.  Consider $S(L)$ as a
  subcomplex of $B_{1/4}$.  We can extend $\pos:S(L)\to S(L)^+$ and
  $\neg:S(L)\to S(L)^-$ to maps $\pos'$ and $\neg'$ defined on
  $B_{1/4}^{(k+1)}\cup S(L)$; since $S(L)^\pm$ is homeomorphic to $L$,
  which is $k$-connected, there is no obstruction to constructing this
  extension and the extensions can be chosen to be Lipschitz.

  Then we can extend $\QQ$ by letting
  $$\QQ(v\cdot x)=\begin{cases} 
    v\cdot s_{|h(v)|}(\neg'(x)) & \text{ if $h(v)>0$} \\
    v\cdot s_{|h(v)|}(\pos'(x)) & \text{ if $h(v)<0$}.
  \end{cases}$$
  for all $x\in B_{1/4}^{(k+1)}$ and $v\in A\setminus H$.
  Again this extension may increase the Lipschitz constant, but
  we still have $\Lip(\QQ|_{h^{-1}([-t,t])})\le c t+c$.
\end{proof}

We will gently abuse notation so that if $\alpha$ is a chain or cycle then we will write the push-forward map $\QQ_\sharp(\alpha)$ as simply $\QQ(\alpha)$,
and similarly for $\P_r$.

\begin{remark}[Pushing and admissible maps]\label{push-admiss}
We note that the pushing maps can be applied to homotopical fillings:
a filling of an admissible (in the sense of \cite{snowflake}) $k$-sphere by an admissible ball
is contained in the $(k+1)$-skeleton, so it can be composed with a
pushing map to get a new filling.  
The volume of the new filling 
is controlled by the Lipschitz constant of the pushing map, and one
can approximate it using an appropriate 
variant of the Deformation Theorem to get a new admissible
filling of the original sphere whose number of cells is controlled.
\end{remark}


\section{Dehn functions in Bestvina-Brady groups}\label{sec:BB}

The kernel $\Hgamma$ of the height map acts geometrically on the zero level set $Z$,
so the Dehn function $\delta_H$ measures the difficulty of filling in $Z$.

When $\Hgamma$ is of type $F_{k+1}$, 
  the results of
  \cite{bestvina-brady} imply that $L$ is $k$-connected.
This means that we have a height-pushing map
$\QQ: Y\cup X^{(k+1)}\to Z$ as defined in the previous section (Lemma~\ref{lem:Q}).
The fact that it is $O(t)$-Lipschitz on the 
part of $X$ up to height $t$ immediately yields bounds on higher Dehn functions;
we will see below that these bounds turn out to be sharp.

\begin{thm}[Dehn function bound for kernels]\label{thm:bestBradUpper}
  If $H=H_\Gamma$ is a Bestvina-Brady group and $H$ is type $F_{k+1}$,
  then $$\delta^{k}_{H}(\V)\preceq \V^{2(k+1)/k}.$$
\end{thm}

The proof is straightforward:  push the CAT(0) filling (Proposition~\ref{prop:wenger})
into the zero level
set, and observe that its volume can't have increased too much while pushing.

\begin{proof}
  Let $a$ be a Lipschitz $k$-cycle in $Z$ of mass at most $\V$, and let
  $t_a=\V^{1/k}$.  By Federer-Fleming approximation 
  (Theorem~\ref{thm:fedflem}), we may assume that $a$ is
  supported in $Z^{(k)}=X^{(k+1)}\cap Z$.    
  
  We know that $\QQ$ has Lipschitz constant $ct+c$ on heights up to $t$.
Since $X$ is CAT(0),  
there is a constant $m>0$ and a chain $b\in
  \CLip_{k+1}(X)$ such that $\mass b\le m t_a^{k+1}$ and $b$ is
  supported in a $mt_a$-neighborhood of $\supp a$.  In particular, the height is bounded:
  $h(\supp b)\subset [-mt_a,mt_a]$.  
  Approximating again, 
  we may assume that $b$ is supported in
  $X^{(k+1)}$.  Then $b'=\QQ(b)$ is a $(k+1)$-chain in $Z$ whose
  boundary is $a$, and 
  $$\delta^k_H(a)\preceq \mass b' \preceq (c m t_a+c)^{k+1} \cdot \mass b
  \preceq t_a^{2k+2}=\V^{\frac{2(k+1)}{k}}.$$\end{proof}

This recovers a theorem of Dison \cite{dison} in the case $k=1$. 

Brady \cite{bradybook} constructed examples of Bestvina-Brady groups
with quartic ($\V^4$) Dehn functions, showing that this upper bound is sharp when $k=1$.  
We next generalize these
examples to find Bestvina-Brady groups with large higher-order Dehn
functions, showing that the upper bound in the previous theorem is sharp for all $k$.

\begin{definition}[Orthoplex groups]
Recall that a $k$-dimensional orthoplex (also known as a cross-polytope) is the 
join of $k+1$ zero-spheres.  The standard $k$-orthoplex is the polytope in $\R^{k+1}$
whose extreme points are $\pm e_i$ for the standard basis vectors $\{e_i\}_{i=0}^k$.
For $k\ge 0$, 
we call a graph $\Gamma$ a $k$-{\em orthoplex graph}, and its associated group $\Agamma$
a $k$-{\em orthoplex group}, if the flag complex $L$ on $\Gamma$ has the following properties:
\begin{itemize}
\item $L$ is a $(k+1)$-complex that is a triangulation of a $(k+1)$-dimensional ball;
\item the boundary of $L$ is isomorphic to a $k$-dimensional orthoplex;
\item there exists a top-dimensional simplex in $L$ whose closure is
  contained in the interior of $L$.  We call this a {\em strictly
    interior} simplex.
\end{itemize}
The boundary is isomorphic as a complex to the standard orthoplex, 
so it has $2(k+1)$ vertices that we label by $a_i,b_i$ for $0\le i\le k$, 
where $a_i$ corresponds to $e_1$ and $b_i$ to $-e_i$.
\end{definition}

For example, a path with at least three edges is a $0$-orthoplex graph.

\begin{figure}[ht]
\includegraphics{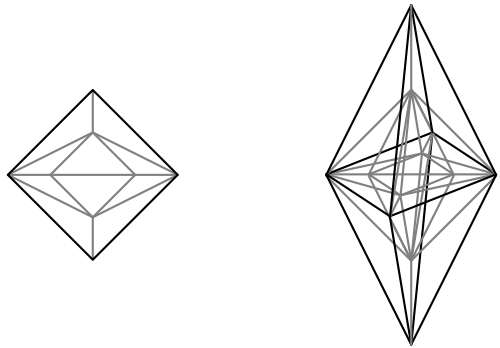}
\caption{The figure on the left is a $1$-orthoplex graph.  It is shown in \cite{bradybook}
that this $H_\Gamma$ has quartic Dehn function.
On the right is a $2$-orthoplex graph.  Note that a copy of the $1$-orthoplex graph
appears on the ``equator'' of the $2$-orthoplex example; similarly, the $1$-orthoplex
graph contains an ``equatorial'' path of length three, which is a $0$-orthoplex graph.  
These are  the simplest symmetric examples of $k$-orthoplex graphs for $k=1,2$.}
\end{figure}

If $A$ is a $k$-orthoplex group, then the flag complex $L$ is a triangulated ball by definition; 
it follows that the associated Bestvina-Brady group $H$ is of finite type~\cite{bestvina-brady}.

\begin{theorem}[Kernels of orthoplex groups have hard-to-fill spheres]\label{dehnsharp}
  If $\Agamma$ is a $k$-orthoplex group, 
  then 
  $$\delta^k_{H}(\V)\succeq \V^{2(k+1)/k}.$$
\end{theorem}

\begin{proof}
We write $A$, $X$, $Z$, and $H$ as usual, with $h:X\to \R$ the height function.
  Let $\gamma_i$, $0\le i\le k$ be the bi-infinite geodesic along edges
  of $X$ such that $\gamma_i(0)=e,$
  $$\gamma_i|_{\R^+}=a_ib_ia_ib_i\dots \quad ; \qquad \gamma_i|_{\R^-}=b_ia_ib_ia_i\dots;$$
  that is, $\gamma_i(-1)=b_i$, $\gamma_i(-2)=b_ia_i$, etc.
  The idea of this proof is that the $\gamma_i$ span a flat that is not collapsed 
  very much by being pushed down to $Z$,
  so that we can get quantitative control on the filling volume in $Z$ for spheres 
  coming from that flat.

In this proof we will adopt the notation that $x=(x_0,\ldots, x_k)\in \Z^{k+1}$.
  Let $F:\Z^{k+1}\to A$ be given by
  $$F(x)=\prod_{i=0}^k \gamma_i(x_i);$$
  note that $\gamma_i(n)$ commutes with $\gamma_j(m)$ for all $i\ne j$
  and all $n,m\in \Z$ and that the image of $F$ lies in the
  non-abelian subgroup $\langle a_i,b_i \rangle$.  The image of $F$ forms the set
  of vertices of a (nonstandard) $(k+1)$-dimensional flat $\bar{F}$.
  This flat $\bar F$ is entirely at non-negative height, with a unique vertex, 
  $F(0, \dots ,0)$, at height zero.
  Since $h(F(x))=\sum |x_i|$, if $r>0$, then the part of $\bar{F}$ at
  height $r$ is an orthoplex and is homeomorphic to $S^k$.
  Define a $k$-sphere $\Sigma_r$ to be a translate of this sphere back to height zero:
  $$\Sigma_r:=[(a_0)^{-r}\bar{F}]\cap Z.$$  We will regard $\Sigma_r$ as a Lipschitz 
  $k$-cycle and show that  it is difficult to fill in $Z.$  
  
  \begin{figure}[ht]
\begin{tikzpicture}[scale=1/2]
\draw[->](-2,0) -- (9,0);
\draw[->] (0,-2) -- (0,9); 
\draw[->](-2,-14/9) -- (9,7);
\draw[->] (-14/9,-2) -- (7,9); 
\draw [very thick, red] (8,0)--(36/9,28/9)--(0,8)--(28/7,36/7)--cycle;
\foreach \a in {0,2,4,6} 
{\node at (\a,-1) [above right] {$b_1$};
\node at (1+\a,-1) [above right] {$a_1$};
\node at (0,\a) [above left] {$a_1$};
\node at (0,1+\a) [above left] {$b_1$};}
\end{tikzpicture}
\includegraphics[height=2in]{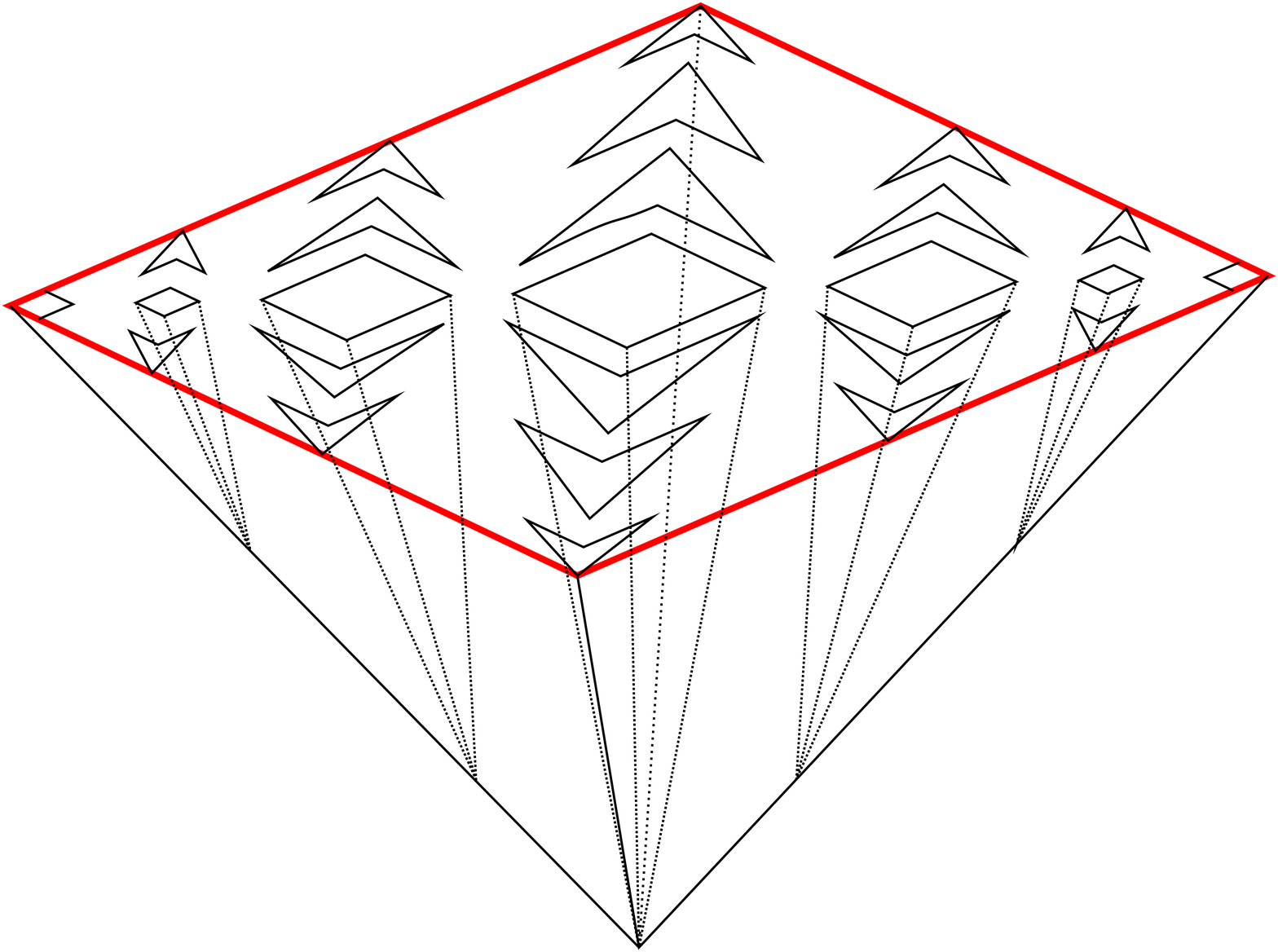}
\caption{The $k=1$ case.
The vertical and horizontal rays in the figure on the left fit together to form the geodesic
$\gamma_1$, and the other two rays fit together to form $\gamma_0$.
The four quadrants glue together to form the plane $\bar F$ (a nonstandard 2-flat).
In red we see the orthoplex $\bar F \cap h^{-1}(r)$ 
whose translate is a hard-to-fill sphere in $Z$.  It has a unique efficient filling in $\bar F$ 
(of area $O(r^2)$), given by coning to the origin.
The figure on the right depicts the pushing of that filling to $Z$,
with signed copies of $L$ shown.  That filling has area $O(r^4)$.
\label{ortholink}}
\end{figure}
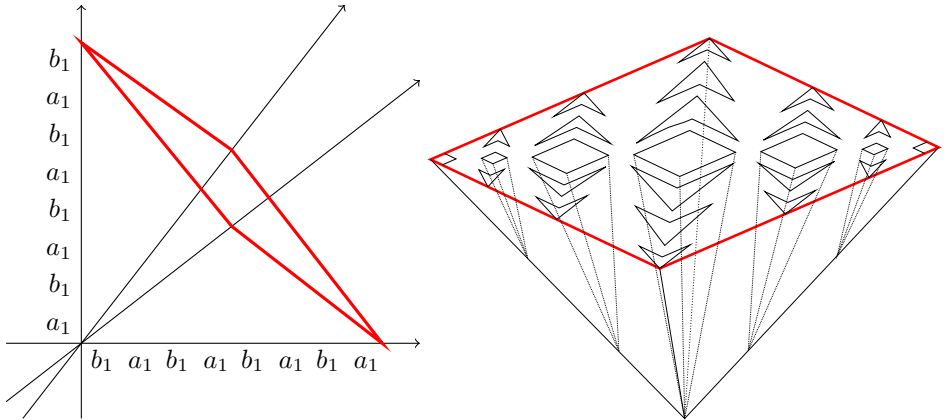

Recall that $Z$ is a $(k+1)$-connected $(k+1)$-dimensional complex.
If $a$ is a cellular $k$-cycle in $Z$, it has a unique cellular
filling $b_0$, and this filling has minimal mass among Lipschitz
fillings.  If $b=\sum b_i \sigma_i$ is a Lipschitz $(k+1)$-chain
filling $a$, where $b_i\in \Z$ and $\sigma_i$ are maps from the
$(k+1)$-simplex to $Z$, then it may have larger mass than $b_0$,
because the different simplices making up $b$ may partially cancel.
This cancellation, however, is the only way that $b$ may differ from
$b_0$.  Since the boundary of $b$ is in the $k$-skeleton of $Z$, the
degree with which it covers any $(k+1)$-cell is well defined, and this
degree must equal the corresponding coefficient of $b_0$.  In
particular, the mass of the parts of $b$ which do not cancel provide a
lower bound for $\|b_0\|$.  Thus, if one of the $\sigma_i$ is disjoint
from all the others then $|b_i|\mass \sigma_i$ is a lower bound on the
mass of any filling.  We will use this general fact to show that
$\Sigma_r$ is hard to fill.

We first use the height-pushing map to find a chain in $Z$ filling $\Sigma_r$;
then we show that this chain is large.  To construct the chain in $Z$ note that 
the cycle $\Sigma_r$ comes from the sphere of radius $r$ in the flat $\bar F$, so it 
has an obvious filling $T_r$  from the ball it bounds in $\bar F$.  Formally, we define
$T_r$  as the $(k+1)$-chain
$$T_r:=(a_0)^{-r}\bar{F}\cap h^{-1}([-r,0]).$$

Before we can apply the height-pushing map, we must perturb $T_r$ so that it 
misses all vertices of $X$ of nonzero height.
In the perturbation $T'_r$, neighborhoods of the vertices of nonzero height 
are replaced with copies of $L$, as follows.  Recall that the link $S(L)$ consists
of signed copies of the simplices of $L$ (Remark~\ref{rem:signedL}); for each simplex of 
$L$ with vertices $v_0,\dots, v_d$, there are $2^{d+1}$ simplices in $S(L)$,
with vertices $\pm v_0,\dots, \pm v_d$.    
For any vertex $v\in T_r$, the link 
$v \cdot S_{1/4} \cap T_r$ is an orthoplex:  it is a join 
of $k+1$ zero-spheres, and the $i$th zero-sphere is labeled by $a_i$ and $b_i$ 
with some signs.  So in $S(L)$, there exists 
a copy of $L$ with this orthoplex as its boundary (in fact there are many, each 
specified by a choice of signs on the interior vertices of $L$).   Perform a surgery at each
vertex $v$ of $T_r$, replacing the $1/4$-neighborhood of $v$ in $T_r$ with a copy 
of $L$ in this way.  The modified filling lives in $Y$, and we call it $T'_r$.

Now,  push the perturbed filling $T'_r$ into $Z$.
The result, $T''_r=\QQ(T'_r)$, is a chain that fills $\Sigma_r$ in
$Z$ (see Figure~\ref{ortholink}).  By the remark above, it suffices to find a lower bound on the
size of the uncanceled pieces of this filling.  

To get such a bound, we only need to consider images under 
$\QQ$ of link faces of $T'_r$, 
since original faces of $Y$ are sent to lower-dimensional pieces.
All of these link faces occur as part of a copy of $L$.  Recall
that $\QQ$ takes each copy of $L$ in its domain to a scaled copy (with some
orientation) of $S(L)^+$ or $S(L)^-$ in  $Z$ (Remark~\ref{rem:signedL}); 
copies with positive
height go to $S(L)^-$ and copies with negative height go to 
$S(L)^+$.  We can thus  write, with $x=(x_0,\ldots,x_k)$, 
\begin{equation}\label{eq:betaSum}
  T''_r=\sum_{j=0}^r  \ \ \sum_{\sum |x_i|=j} (-1)^j(a_0)^{-r}F(x) s_{r-j}(\lambda),
\end{equation}
where $\lambda$ is the fundamental class of $S(L)^+$.
We are trying to estimate $\|T''_r\|$, but some of the terms in
\eqref{eq:betaSum} may cancel.  We will obtain a lower bound on
$\|T''_r\|$ by showing that many of the scaled simplices making up
the sum are disjoint from all other cells of $T''_r$.

  First, note that if $\sigma$ and $\sigma'$ are two different
  $k$-simplices of $S(L)^+$, and $g,g'\in A$ and $t,t'>0$, then 
  $g s_t(\sigma)$ and $g' s_{t'}(\sigma')$ intersect in at most a
  $(k-1)$-dimensional set.  We thus only need to consider the case of
  overlap between scaled copies $g s_{t}(\sigma)$ and $g' s_{t'}(\sigma)$ of
  the same $\sigma$.  Let $\sigma$ be a strictly interior $(k+1)$-simplex in 
  $L$, with vertex set $S=\{g_0,\dots, g_k\}$; note that $S$ is disjoint from 
  $\{a_1,\dots, a_k, b_1,\dots, b_k\}$.  Let $A_S=\langle g_0,\dots, g_{k}\rangle$.  We claim that the
  scaled copies of $\sigma$ in \eqref{eq:betaSum} are disjoint, and so
  none of them is canceled in the sum.  If $v$ and $v'$ are vertices
  of scaled copies of $\sigma$, based at $x\in \bar{F}$ and $x'\in
  \bar{F}$ respectively, then
  $$v=(a_0)^{-r}x \prod_i g_i^{t_i} \quad ; \qquad
  v'=(a_0)^{-r}x' \prod_i g_i^{t'_i}$$
  for some vertex $x$ of $\bar{F}$ and some $t_i\in \Z^+$.  Note that
  $x\in \langle a_i,b_i\rangle$ and $\prod_i g_i^{t_i'}\in
  A_S$.  Since $A_S\cap \langle a_1,\dots, a_k, b_1,\dots, b_k \rangle=\{0\}$,
  if $v=v'$, then $x=x'$, so no two distinct scaled copies of $\sigma$ intersect.

  Consequently, these terms do not cancel in $T''_r$, and as we mentioned
  previously, the mass of any filling of $\Sigma_r$ is bounded below by the mass
  of these terms.  Thus
$$    \delta^k_X(\Sigma_r)  \ge \sum_{j=0}^r  \   \sum_{\sum |x_i|=j} \mass(s_{r-j}(\sigma))  
     \ge \sum_{j=0}^r  \    \sum_{\sum |x_i|=j} \mass(\sigma)(r-j)^{k+1}    \ge c r^{2k+2} $$
  for some $c>0$ depending only on $k$.  Since $\mass(\Sigma_r)\asymp r^{k}$ 
  and $r$ was arbitrary, we have $\delta^k_H(\V)\succeq \V^{\frac{2k+2}{k}}$ as desired.
  \end{proof}

\section{Higher divergence in RAAGs}\label{sec:higherdiv}

There are still many subgroups of right-angled Artin groups for which the Dehn
function is unknown, and the higher divergence functions have a
similar level of difficulty.  For the groups themselves, however, we can get sharp 
bounds on the possible rates of divergence.

\begin{theorem}[Higher divergence in RAAGs]\label{divk}  
For $0\le k\le \divdim(A_\Gamma)$,
  $$  r^{k+1} \preceq  \div^k(\Agamma) \preceq r^{2k+2}.$$
  The upper and lower bounds are sharp:  for every $k$ there are examples of
  right-angled Artin groups realizing these bounds.
\end{theorem}

In the next section, we will give a sharper result in the case $k=0$.

In \S\ref{sec:euclidean example} we saw that the general lower bounds are realized 
by free abelian groups.  We divide the rest of the theorem into several pieces:  the 
general lower bounds, the general upper bounds, and a construction of groups whose 
divergence realizes the general upper bounds.

\begin{proposition}[RAAG lower bounds]
If $k\le \divdim(A_\Gamma)$,
  then $  r^{k+1} \preceq  \div^k(\Agamma)$.
\end{proposition}

\begin{proof}
The Dehn function of $\Agamma$, evaluated at $\alpha r^k$, is a lower
bound for $\ddiv_\rho^k(\alpha r^k, r)$.  This is because there are cycles in $\Agamma$
of mass at most $\alpha r^k$ whose most efficient fillings have mass arbitrarily close to
$\delta^k_{A}(\alpha r^k)$, and these can be translated to be $r$-avoidant.  
Since $\Gamma$ has a clique of $k+1$ vertices, $\Agamma$ retracts onto
a subgroup $\Z^{k+1}\subset \Agamma$.  Thus
$$\delta^k_{\Agamma}(\V)\succeq \delta^k_{\Z^{k+1}}(\V)\succeq \V^{\frac{k+1}{k}}.$$
\end{proof}

\begin{proposition}[RAAG upper bounds]\label{genupperbounds}
If $k\le \divdim(A_\Gamma)$, then 
  $\div^k(\Agamma)\preceq r^{2k+2}$.
\end{proposition}

\begin{proof}
  By Remark~\ref{compare}(1), it suffices to show that there is a $c$ such that for 
  sufficiently large $r$,
  $$\ddiv_1^k(\V, r)\le c \V^{\frac{k+1}{k}} r^{k+1}.$$
  Let $a\in \Zcell_k(X_\Gamma)$ be an $r$-avoidant $k$-cycle and let
  $\V=\|a\|$.  Since $X_\Gamma$ is CAT(0), there is a $(k+1)$-chain
  $b\in \Ccell_{k+1}(X_\Gamma)$ such that $\partial b=a$ and $\|b\|_1\preceq
  \V^{\frac{k+1}{k}}$.  We will consider $b$ as a Lipschitz chain and
  push it out of $B_r$.
  
  Let $\P_r$ be the map constructed in Lemma~\ref{lem:P}; this map is
  $O(r)$-Lipschitz.  The image $\P_r(b)$ is an
  $r$-avoidant filling of $a$, and there is a $c$ such that 
  $$\mass (\P_r(b))\le (\Lip \P_r)^{k+1} \mass b\le
  cr^{k+1} \V^{\frac{k+1}{k}},$$
  as desired.
\end{proof}

The final step of Theorem~\ref{divk} is to construct groups that have
the stated divergences.  The key to this construction is to use the
connection between divergence functions and the Dehn functions of
Bestvina-Brady groups; large portions of $S_r$ can be embedded in the level sets corresponding to
Bestvina-Brady groups, so avoidant fillings can be converted
into fillings in Bestvina-Brady groups.

\begin{theorem}[Sharpness of upper bounds]\label{divsharp}
  If $\Agamma$ is a $k$-orthoplex group, then $\div^k(\Agamma)\asymp r^{2k+2}$.
\end{theorem}

\begin{proof}
  By Remark~\ref{compare}(1), it suffices to show that there is a $c_k>0$ such
  that for all $0<\rho \le 1$,
  $$\ddiv_\rho^k(c_k r^k, r)\succeq r^{2k+2}.$$

  Recall that when $A_\Gamma$ is a $(k+1)$-dimensional orthoplex group,
  we constructed cycles in $Z_\Gamma$  by defining a flat $\bar{F}$ and
  considering the intersections
  $$\Sigma_r:=((a_0)^{-r}\bar{F})\cap Z.$$
  These have mass $c_k r^k$ for some $c_k>0$ and require 
 mass $\asymp r^{2k+2}$ to fill.  Let 
  $$\Sigma'_{r}:=(a_0)^{r}\Sigma_r=\bar{F}\cap h^{-1}(r).$$
  This is an $r$-avoidant cycle of volume $c_k r^k$ and we will show
  that every $\rho r$-avoidant filling of $\Sigma'_{r}$ has volume
  $\succeq r^{2k+2}$.

  We first define a retraction $\pi_t:X_\Gamma\setminus B_{t}\to S_{t}$
  for all $t$.  If $x\in X_\Gamma\setminus B_t$, there is a unique
  CAT(0) geodesic  from
  $x$ to $e$.  Let $\pi_t(x)$ be the intersection of $S_{t}$ with this
  geodesic.  (As before, we will also 
  write $\pi_t$ for the induced map $(\pi_t)_\sharp$ on chains and cycles.)
  This is clearly the identity on $S_t$, and since $B_t$ is
  convex, it is 1-Lipschitz (distance-nonincreasing).  
  Furthermore, if $x\in \bar{F}$,
  then since $\bar{F}$ is a flat, the geodesic from $x$ to $e$ is a
  straight line in $\bar{F}$, and
  $\pi_t(\Sigma'_{r})=\pi_t(\Sigma'_{t})$ for all $t\le  r$.

  Consider the map $\abs \circ \pi_t$.  We claim that this is a
  1-Lipschitz retraction from $X_\Gamma\setminus B_t$ to
  $S_t\cap h^{-1}(t)$.  If $x\in S_t\cap h^{-1}(t)$, then
  $\pi_{t}(x)=x$, and $\abs(x)=x$, so this map is a retraction, and
  since $\abs$ and $\pi_t$ are each 1-Lipschitz, the composition is as well.  
  For all $t\le r$ we have 
  $$\abs \circ \pi_t(\Sigma'_r)=\Sigma'_t.$$

  Fix an arbitrary $0<\rho\le 1$ and 
  let $b$ be a $\rho r$-avoidant $(k+1)$-chain whose boundary is
  $\Sigma'_r$.  Then $b'=\abs\circ \pi_{\rho r}(b)$ is a
  chain in $S_{\rho r}\cap h^{-1}(\rho r)$ whose boundary is $\partial
  b'=\Sigma'_{\rho r}$.  Its translate $(a_0)^{-\rho r}b'$ is a chain
  in $Z$ whose boundary is $\Sigma_{\rho r}$, so 
  $$\mass b'\ge \delta^k_H(\Sigma_{\rho r})\succeq r^{2k+2}.$$
  Since $\abs \circ \pi_t$ is 1-Lipschitz, $\mass b\ge \mass
  b'$, so
  $$\ddiv_\rho^k(c_k r^k, r)\succeq r^{2k+2}$$
  as desired.
\end{proof}

\section{A refined result for divergence of geodesics}\label{sec:div0}

The case $k=0$ gives a quantitative measure of how fast geodesics spread apart.
Here, the answer is completely determined by whether or not the group
is a direct product.  Note that $\Agamma$ is a direct product if and only if the vertices 
can be partitioned into two nonempty subsets $A$ and $B$ such that each vertex of
$A$ is joined to each vertex of $B$ by an edge of $\Gamma$.  Equivalently, $A_\Gamma$
is a direct product if and only if the complement of $\Gamma$ is not connected.

The following result also appears in \cite{behrstock-charney}, with a completely 
different proof.

\begin{theorem}[Divergence of geodesics in RAAGs]\label{div0raag} 
For a right-angled Artin group $\Agamma$, $\div^0$ exists if and only if the defining graph 
$\Gamma$ is connected.
In this case, $\div^0(\Agamma)\asymp r$ if and only if $\Agamma$ is a nontrivial 
direct product, and $\div^0(\Agamma)\asymp r^2$ otherwise.
\end{theorem}

Note that if $\Gamma$ is not connected, then $\Agamma$ has infinitely many ends, so $\div^0$ is
not defined.  
We have already established (Theorem~\ref{divk} with $k=0$)
that $r\preceq \div^0(\Agamma) \preceq r^2$.  
We proceed by considering the presence of a product structure.

\begin{lemma}[Linear if direct product]\label{lem:LinearUpper}
If $\Agamma=H \times K$ is a direct product of nontrivial factors, then $\div^0(\Agamma) \asymp r$.  
\end{lemma}

\begin{proof} 
Writing elements of $\Agamma$ as ordered pairs, let $(h_1, k_1)$ and $(h_2, k_2)$ be elements of $\Agamma$ 
of length $r$, that is, $|h_i|_H + |k_i|_K =r$.  There exists a $u \in H$   such that $|u|_H \leq r$ and  
$|h_1u|_H = r$.  Similarly there exists an element 
$v \in K$ such that $|v|_K \leq r$ and  $|k_2v|_K \geq r$.  
Now
the vertices representing 
the elements 
\[
(h_1, k_1), \,
(h_1 u, k_1), \,
(h_1 u, e), \, 
(h_1 u, k_2 v), \,
(e, k_2 v), \,
(h_2, k_2 v), \,
(h_2, k_2)
\]
lie on or outside 
$B_r$, and successive elements of the sequence can be connected by $r$-avoidant paths, each of which has length at most $r$.  Thus any two vertices on 
$S_r$ can be connected by a $r$-avoidant path of length at most $6r$, and the lemma follows. 
\end{proof}

In fact the proof uses little about RAAGs; the same result holds for all direct products $H\times K$
where $H$ and $K$ each have the property that every point lies on a geodesic
ray based at $e$.

\begin{lemma}[Quadratic if not direct product]\label{lem:quad-lowe}
If $\Agamma$ is not a direct product, then 
$\div^0(\Agamma) \asymp r^2$.
\end{lemma}

\begin{proof}  We only need to show that $\div^0\succeq r^2$.  
As $\Agamma$ is not a direct product, the complement $\Gamma^c$ of $\Gamma$
is connected.
Choose a closed path in $\Gamma^c$ that visits each vertex (possibly with repetitions)
and let $w=a_1 a_2 \dots a_n$ be the word made up of the generators 
encountered along this path. 
Introduce the symbol $a_{n+1}$ as another name for $a_1$.
Note that for $1 \leq i \leq n$, the vertices $a_i$ and $a_{i+1}$ are not connected by 
an edge in $\Gamma$, 
so the corresponding generators do not commute. 
As a consequence, $w$ is a geodesic word, as is any nontrivial power $w^k$; let 
$\eta$ be the unique
bi-infinite geodesic going through $e$ and all of these powers, so that the letters of $\eta$ 
cycle through the $a_i$.  (Note this is a geodesic with respect to either the word metric 
on $\Agamma$ or the CAT(0) metric on $X$.)  

To prove the lemma, it is enough to show that for each $r$, and each $0< \rho\leq 1$,  any $\rho r$-avoidant path connecting 
$\eta(r)$ and $\eta(-r)$  has length at least on the order of $r^2$.  

Let $\gamma$ be a $\rho r$-avoidant path from $\eta(r)$ to $\eta(-r)$ in $X^{(1)}$.  
Let $\eta_0$ be the segment of $\eta$ between the same two endpoints. 
The concatenation of $\gamma$ and $\eta_0$ taken in reverse is a loop in $X$ labeled by a word representing the identity in 
$\Agamma$ (see Figure~\ref{fig:cors}).
So there exists a van Kampen diagram $\Delta$ whose boundary cycle $\partial \Delta$ is labeled by this word, 
and a combinatorial map $\Delta \to X$ such that $\partial \Delta$ maps to $\gamma \cup \eta_0$.

If $\sigma  \in \partial \Delta$ is an edge mapping to $\eta_0 \cap B_{\rho r}$, then $\sigma$ is in the boundary of 
a $2$-cell of $\Delta$.   In the informal ``lollipop'' language used to talk about these diagrams, $\sigma$ is part of 
the candy.   To see this note that the only way an edge in $\eta_0$ can come from the 
``sticks'' of the lollipop is if it coincides with some edge of $\gamma$, but $\gamma$ was chosen to be $\rho r$-avoidant.

If $y$ is the label on $\sigma $, then $\sigma $ is one end of a $y$-corridor whose
other end is an 
edge of $\partial \Delta$ 
with label $y$, but orientation opposite to that of $\sigma$.  Since the map from 
$\Delta$ to $X$ preserves orientations of 
edge labels, and all the edges of $\eta$ have the same orientation, 
the other end of the corridor must be an 
edge of $\gamma$. 
Thus each edge in $\eta_0\cap B_{\rho r}$ bounds a corridor whose other end is an edge of 
$\gamma$ (see Figure~\ref{fig:cors}).
Since $a_i$ does not commute with $a_{i+1}$  for any $i$, no two of these corridors intersect.

\begin{figure}[ht]
\begin{tikzpicture}[scale=1]
\filldraw (-5,-1) circle (0.05) node [below] {$\eta(-r)$};
\filldraw (6,-1) circle (0.05) node [below] {$\eta(r)$};
\filldraw (0.5,-1) circle (0.05) node [below] {$e$};
\node at (-2.5,.5) {$\cdots$};
\node at (-2.5,4) {$\gamma$};
\node at  (-3.2,-1) [below] {$a_1$};
\node at  (-2.8,-1) [below] {$a_2$};
\node at  (-2,-1) [below] {$a_n$};
\node at  (1.4,-1) [below] {$a_1$};
\node at  (2.2,-1) [below] {$a_i$};
\node at  (2.8,-1) [below] {$a_n$};
\node at  (3.8,-1) [below] {$a_i$};
\node at  (5,-1) [below] {$a_1$};
\node at (4.7,2.9) {$x$};
\begin{scope}[inner sep=0]
\draw [clip] (-5,-1) .. controls (-4,4) and (-3,3) .. (-2,4) .. controls (-1,5) and (1,4.3) .. (2,3.8) 
.. controls  (3,3.5) and (5,4) .. (6,-1) -- cycle;
\draw (-3.4,-1)--(-5.5,3) (-5.1,3)--(-3,-1)--(-4,3) (-2.6,-1)--(-3.6,3);
\draw (-2.2,-1)--(-1.5,5) (-1.8,-1)--(-1.1,5);
\draw (1.2,-1)--(-1,5) (1.6,-1)--(-.6,5) (2,-1)--(1,4.5) (2.4,-1)--(1.4,4.5)  
(2.6,-1)-- node [pos=.5] (M1) {} node [pos=.6] (M2) {} 
(3.5,4) (3,-1)-- node [pos=.5] (N1) {} node [pos=.6] (N2) {} (3.9,4)  
(3.6,-1)--(6,2) (4,-1)--(6.4,2) (4.8,-1)--(6,1) (5.2,-1)--(6.4,1);
\draw (2,4)--(M1)--(N1)--(5,3) (2.4,4)--(M2)--(N2)--(4.6,3)  ;
\end{scope}
\end{tikzpicture}
  \caption{Corridors in $\Delta$.  For each of the corridors from $\eta$ to $\gamma$ (vertical in this picture),
  every cell has a corridor that crosses it with both ends on $\gamma$.
    \label{fig:cors}}
\end{figure}
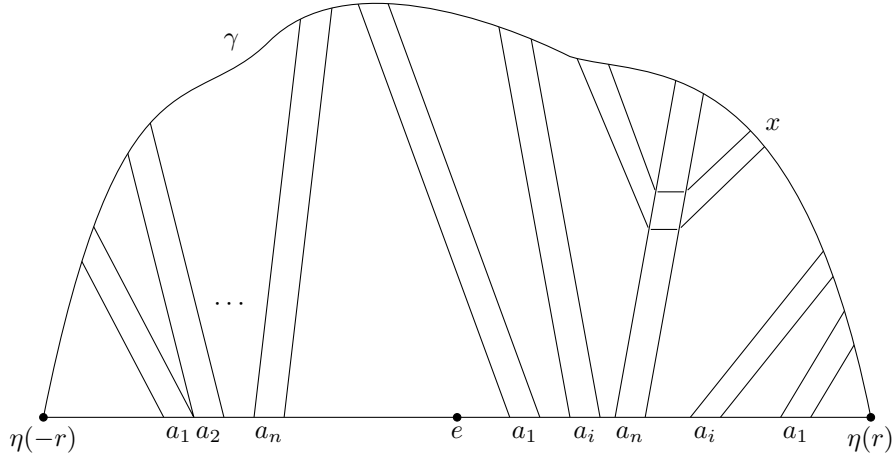

Any $y$-corridor has boundary label
$yvy^{-1}v^{-1}$ for some word $v$, which we call the lateral boundary word of the corridor.  
We now show that by performing some surgeries on $\Delta$ which do not change 
$\partial \Delta$, we may assume that the lateral boundary words of the $a_i$-corridors in $\Delta$ emanating from $\eta$ are 
geodesic words
(i.e., minimal representatives of the group elements that they represent).  

Suppose the boundary word of an $a_i$-corridor is not geodesic. 
Since every word can be reduced to a geodesic by shuffling neighboring commuting pairs 
(see~\cite{herm-meier}),
there has to be a sub-segment of the form $xux^{-1}$ where $x$ is a generator (or its inverse), $u$ is a word, and 
$x$ commutes with the individual letters of $u$.  
Then we can perform the \emph{tennis-ball move} shown in Figure~\ref{tennisball}.  

\begin{figure}[ht]
\begin{tikzpicture}[scale=.6]
\draw [fill=blue!30, dashed] (4,2.3) rectangle (2,4.3);
\draw [dashed] (-2,0)--(2,2.3); \draw (0,0)--(4,2.3)--(4,4.3)--(0,2)--cycle;
\draw [fill=blue!30] (4,4.3)--(2,4.3)--(-2,2)--(0,2)--cycle; \draw [fill=blue!30] (0,0) rectangle (-2,2);
\draw [fill=blue!30] (6.8,-1) rectangle (8,5); \draw [->,very thick] (8.5,2)--(9,2);
\draw (9.4,0) rectangle (13,4); \draw (10.6,0) -- (10.6,4);  \draw (11.8,0) -- (11.8,4);
\end{tikzpicture}
\caption{Here, the long edges are labeled by a word $u$ and the square faces are $x,a_i$-commutators.
In a {\em tennis-ball move},  the shaded disk is replaced with the unshaded disk, noting that the boundary words 
($xux^{-1}a_i x u^{-1} x^{-1} a_i^{-1}$)
are equal.}\label{tennisball}
\end{figure}
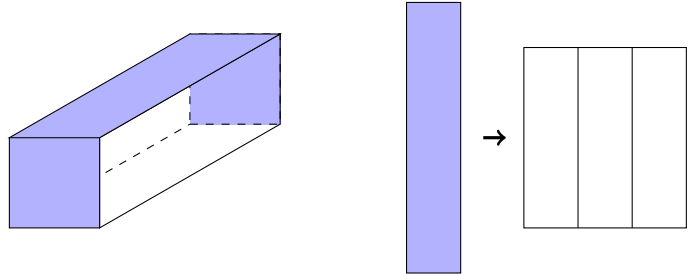

Note that $\partial \Delta$ remains unchanged at the end of such a move, and that the length of the $a_i$-corridor is reduced.  
After performing enough of these moves 
we get a van Kampen diagram, which we also call $\Delta$, in which lateral boundaries of corridors emanating from $\eta$ are geodesics.
There is a map from $\Delta$ to $X$, which agrees with the original one on $\partial \Delta$. 

We now restrict our attention to $a_1$-corridors emanating from $\eta$.  
A $2$-cell in such a corridor has boundary label 
$a_1 x a_1^{-1} x^{-1}$ for some $x$, and is therefore part of an $x$-corridor or annulus.  We claim that this is in 
fact an $x$-corridor which intersects the $a_1$-corridor in exactly
one $2$-cell.  (In particular it is not an annulus.)  
If the intersection contains more than one $2$-cell, then $\Delta$ contains the picture shown in Figure~\ref{fig:intersect}, where $u$ is a geodesic word 
and $v$ is a word (not 
necessarily geodesic) whose individual letters commute with $x$.
The words $u$ and $v$ represent the same group element, and so each generator appearing in the geodesic word $u$ also appears in $v$.  (This follows from~\cite{herm-meier}.) Thus 
 $x$ commutes with the individual letters in $u$, and 
hence with $u$.  This contradicts the fact that the boundary word of the $a_1$-corridor was geodesic, 
and proves the claim. 

\begin{figure}[ht]
\begin{tikzpicture}[scale=.4]
\draw (-2,-5)--(-2,5);  \node at (-1.5,3) [above] {$a_1$};  \node at (-2,0) [left] {$u$};
\draw (-1,-5)--(-1,5);  \node at (-1.5,-3) [below] {$a_1$};  \node at (-1,0) [right] {$u$};
\draw (-5,-3)--(0,-3);  \node at (-2,-2.5) [left] {$x$};
\draw (-5,-2)--(0,-2) (-5,2)--(0,2)  (-5,3)--(0,3);  \node at (-2,2.5) [left] {$x$};
\draw (0,-3) arc (-90:90:3);  \node at (3,0) [right] {$v$};
\draw (0,-2) arc (-90:90:2);  \node at (2,0) [left] {$v$};
\end{tikzpicture}
  \caption{Each crossing corridor intersects the vertical corridor only once.}
  \label{fig:intersect}
\end{figure}
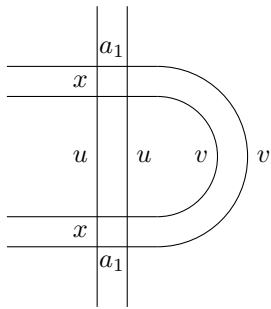

Now $x=a_i$ for some $i$, and since two $a_i$-corridors can't intersect in a $2$-cell, the $x$-corridor is trapped in the region bounded by two $a_i$-corridors 
(see Figure~\ref{fig:cors}), and its two ends are edges of $\gamma$.  
In particular, any two corridors that intersect $a_1$-corridors
emanating from $\eta$ are distinct.  
%
Thus each $2$-cell along each of the $a_1$-corridors emanating from 
$\eta_0$ is part of a corridor whose ends are edges of $\gamma$ and no two of these edges coincide.  Thus the total 
area of the $a_1$-corridors is a lower bound for the length of $\gamma$. 
To obtain a lower bound on this area, note that for each $0\leq j \leq \lfloor r/n\rfloor$, 
there is an $a_1$-corridor whose 
$\eta$-end has vertices whose distances from the origin are $jn$ and $jn+1$ respectively.  Since the other end of 
this corridor is an edge along 
$\gamma$, and $\gamma$ is $\rho r$-avoidant, 
the length of the corridor is at least $\rho r-jn$.  So the total 
area of the $a_1$-corridors is at least 
$\sum_{j=0}^{\lfloor r/n \rfloor} (\rho r - jn) $,
which is on the order of $r^2$. 
\end{proof}

\bibliographystyle{plain}
\bibliography{divk}

\end{document}